\documentclass[11pt]{article}
\usepackage{amsmath,amssymb,natbib}
\usepackage{setspace}
\usepackage{graphicx,color}

\definecolor{DarkBlue}{rgb}{0.1,0.1,0.5}
\definecolor{Red}{rgb}{0.9,0.0,0.1}
\definecolor{purple}{rgb}{0.8,0.0,0.4}

\newtheorem{theorem}{Theorem}

\newtheorem{lemma}{Lemma}

\newtheorem{remark}{Remark}

\def\tr{\mbox{tr}}
\def\bs{\boldsymbol}

\begin{document}
\bibliographystyle{plainnat}

\begin{center}
{\Large{\bf Minimax bounds for sparse PCA with noisy high-dimensional data}}

\vskip.1in Aharon Birnbaum$^*$, Iain M. Johnstone$^\dag$, Boaz Nadler$^\ddag$
\textit{and} Debashis Paul$^\S$

\vskip.1in\textit{\small $*$ Hebrew University of Jerusalem; $\dag$ Stanford
University; $\ddag$ Weizmann Institute of Science; $\S$ University of
California, Davis}

\end{center}

\begin{abstract}
We study the problem of estimating the leading eigenvectors of a
high-dimensional population covariance matrix based on independent Gaussian
observations. We establish a lower bound on the minimax risk of estimators
under the $l_2$ loss, in the joint limit as dimension and sample size increase
to infinity, under various models of sparsity for the population eigenvectors.
The lower bound on the risk points to the existence of different regimes of
sparsity of the eigenvectors. We also propose a new method for estimating the
eigenvectors by a two-stage coordinate selection scheme.
\end{abstract}



\section{Introduction}

Principal components analysis (PCA) is a widely used technique in
reducing dimensionality of multivariate data. A traditional setting
where PCA is applicable involves repeated observations from a
multivariate normal distribution.
Two key theoretical questions are: {\em i) what is the relation between the
sample eigenvectors and the population ones ? and ii) how well can population
eigenvectors be estimated under various sparsity assumptions ?} When the
dimension $N$ of the observations is fixed and the sample size $n$ increases to
infinity, the asymptotic properties of the sample eigenvalues and eigenvectors
are well-known \citep{Anderson1963,Muirhead1982}. Most of this asymptotic
analysis is
based on the fact that the sample covariance approximates well the
population covariance when the sample size is large. However, it is
increasingly common to encounter statistical problems where the
dimensionality of the observations is of the same order of magnitude
as (or even bigger than) the sample size. In such cases, the
sample covariance matrix, in general, is not a reliable
estimate of the population covariance matrix.

To overcome this curse of dimensionality, several works studied the estimation
of the population covariance matrix, under various models of sparsity. These
include the development of banding and thresholding schemes
\cite{BickelL2008a,BickelL2008b,ElKaroui2008,RothmanLZ2009,CaiL2011}, and
analysis of their rate of convergence in the spectral norm. More recent works,
such as \cite{CaiZZ2010} and \cite{CaiZ2011} established the minimax rate of
convergence under the matrix $l_1$ norm and the spectral norm, and its
dependence on the assumed sparsity level.


In contrast to these works, that studied estimation of the population
covariance matrix, in this paper we consider a related but different problem,
namely, the estimation of its leading eigenvectors. The interest in comparing
these two problems is partially due to the fact that, when the population
covariance is a low rank perturbation of the identity, which is a primary focus
of this paper, sparsity of the eigenvectors corresponding to the non-unit
eigenvalues implies sparsity of the whole covariance. Note that consistency of
an estimator of the whole covariance matrix also implies convergence of its
leading eigenvalues to their population counterparts. If the gaps between the
neighboring distinct eigenvalues remain bounded away from zero, it also implies
convergence of the corresponding eigen-subspaces \cite{ElKaroui2008}. Moreover,
for population eigenvalues with multiplicity one and gaps with neighboring
eigenvalues bounded away from zero, the upper bounds for the whole covariance
estimation under the spectral norm, derived in \cite{BickelL2008b} and
\cite{CaiZ2011}, also yield an upper bound on the rate of convergence of the
corresponding eigenvectors under the $l_2$ loss. These works, however, did not
study the following fundamental problem, considered in this paper: \textit{How
well can the leading eigenvectors be estimated, namely, what are the minimax
rates for eigenvector estimation ?}


We formulate this eigenvector estimation problem under the well-studied
``spiked population model'' which assumes that
\begin{itemize}
\item[(*)] the eigenvalues of the population covariance matrix $\Sigma$ are
$$
\lambda_1 + \sigma^2, \ldots, \lambda_M + \sigma^2, \sigma^2,
\ldots,\sigma^2,
$$
for some $M \geq 1$, where $\sigma^2 > 0$ and $\lambda_1 >\lambda_2 > \cdots >
\lambda_M > 0$.
\end{itemize}
This is a standard model in several scientific fields, including for example
array signal processing (e.g. see \cite{vanTrees2002}) where the observations
are modeled as the sum of an $M$-dimensional random signal and an independent,
isotropic noise. It also arises as a latent variable model for multivariate
data, for example in factor analysis \citep{Jolliffe2002,TippingB1998}. The
assumption that the leading $M$ eigenvalues are distinct is made to simplify
the analysis, as it ensures that the corresponding eigenvectors are
identifiable up to a sign change.
The assumption that all remaining eigenvalues are equal is not crucial as our
analysis can be generalized to the case when these are only bounded by
$\sigma^2$.
Asymptotic properties of the eigenvalues and eigenvectors of the sample
covariance matrix under this model, in the setting when $N/n \to c \in
(0,\infty)$ as $n \to \infty$, have been studied by \cite{Baik2006},
\cite{Nadler2008}, \cite{Onatski2006} and \cite{Paul2007}, among others. A
conclusion of these studies is that when $N/n \to c > 0$, the eigenvectors of
standard PCA are inconsistent estimators of the population eigenvectors.

In analogy to the sparse covariance estimation setting, several works
considered various models of sparsity for the leading eigenvectors and
developed improved sparse estimators. For example \cite{WittenT2009} and
\cite{ZouHT2006}, among others, imposed $l_1$-type sparsity constraints
directly on the eigenvector estimates and proposed optimization procedures for
obtaining them. \cite{ShenH2008} suggested a regularized low rank approach to
sparse PCA. The consistency of the resulting leading eigenvectors was recently
proven in \cite{ShenSM2011}, using a formulation of sparsity in which the
sample size $n$ is fixed while $N \to \infty$. \cite{dAspremontEJL2008}
suggested a semi-definite programming (SDP) problem as a relaxation to the
$l_0$-penalty for sparse $\Sigma$. Assuming a single spike, \cite{AminiW2008}
studied the asymptotic properties of the leading eigenvector of the covariance
estimator obtained by \cite{dAspremontEJL2008}, in the joint limit as both
sample size and dimension tend to infinity.
Specifically, \cite{AminiW2008} considered a leading eigenvector with exactly
$k \ll N$ nonzero entries all of the form $\{-1/\sqrt{k},1/\sqrt{k}\}$. For
this hardest subproblem in the $k$-sparse $l_0$-ball,  \cite{AminiW2008}  first
derived information theoretic lower bounds, and then, under the assumption that
the SDP problem has a rank one solution, proved that it attains the optimal
rate of convergence.



In this paper, in contrast, following \cite{JohnstoneL2009} we study the
estimation of the leading eigenvectors of \(\Sigma\) assuming that these are
approximately sparse, with a bounded $l_q$ norm. Under this model,
\cite{JohnstoneL2009} developed an estimation procedure based on coordinate
selection by thresholding the diagonal of the sample covariance matrix,
followed by the spectral decomposition of the submatrix corresponding to the
selected coordinates. \cite{JohnstoneL2009} further proved consistency of this
estimator assuming dimension grows at most polynomially with sample size, but
did not study its convergence rate. Since this estimation procedure is
considerably simpler to implement and computationally much faster than the
$l_1$ penalization procedures cited above, it is of interest to understand its
theoretical properties. More recently, \cite{Ma2011} developed a related scheme
named ITSPCA (iterative thresholding sparse PCA) which is based on repeated
application of filtering, thresholding and orthogonalization steps that result
in sparse estimators of the subspaces spanned by the leading eigenvectors.
He also proved consistency and derived rates
of convergence of the proposed estimator under appropriate loss functions and
sparsity assumptions.

In this paper, which is partly  based on the Ph.D. thesis \cite{Paul2005} and
\cite{PaulJ2007}, we study the estimation of the leading eigenvectors of
$\Sigma$ within the framework of \cite{JohnstoneL2009}, but with an arbitrary
number of spikes (i.e., $M \geq 1$) whose corresponding eigenvectors all belong
to appropriate $l_q$ spaces.
Our analysis thus extends the setting studied in \cite{JohnstoneL2009} and
complements the work of \cite{AminiW2008} that considered the $l_0$-sparsity
setting. For simplicity,
we assume Gaussian observations in our analysis. However, up to multiplicative
constants, the bounds on the minimax rate reported in this paper continue to
hold under a relaxed assumption of sub-Gaussian tail behavior for the
probability distributions of the random variables.

The main contributions of this paper are as follows. First, we establish lower
bounds on the rate of convergence of the minimax risk for any eigenvector
estimator under the $l_2$ loss. This analysis points to three different regimes
of sparsity, which we denote as \textit{dense, sparse, and ultra-sparse}, each
having a different rate of convergence. We show that in the ``dense'' setting
(as defined in Section \ref{sec:lower_bound}), the standard PCA estimator
attains the optimal rate of convergence, whereas in sparse settings it is not
even consistent.  Next, we show that while the diagonal thresholding scheme of
\cite{JohnstoneL2009} is consistent under these sparsity assumptions, in
general, it is not rate optimal. This motivates us to propose a new method
(Augmented Sparse PCA, or ASPCA) for estimating the eigenvectors that is based
on a two-stage coordinate selection scheme, and is a refinement of the
thresholding scheme of \cite{JohnstoneL2009}. While beyond the scope of this
paper, it is possible to show that in the ultra-sparse setting, both our ASPCA
procedure, as well as the method of \cite{Ma2011}   achieve the lower bound on
the minimax risk obtained by us, and are thus  rate-optimal procedures. There
is an intermediate region where a gap exists between the current lower bound
and the upper bound on the risk. It is an open question whether the lower bound
can be improved in this scenario, or a better estimator can be derived. Table
\ref{T:comparison} provides a comparison of the lower bounds and rates of
convergence of  various estimators.

The theoretical results also show that under comparable scenarios, the optimal
rate of convergence for eigenvector estimation,  \(O((\log N/n)^{-(1-q/2)})\)
(under squared-error loss) is faster than the optimal rate for covariance
estimation, \(O((\log N/n)^{-(1-q)})\) (under squared operator norm loss), as
obtained by \citep{BickelL2008b} and \cite{CaiZ2011}. Finally, we emphasize
that to obtain good finite-sample performance for both our two-stage scheme, as
well as for other thresholding methods, the exact thresholds need to be
carefully tuned. This issue and the detailed theoretical analysis of the ASPCA
estimator is beyond the scope of this paper, and will be presented in a future
publication. After this paper was completed, we learned of \cite{VuL2012},
which cites \cite{PaulJ2007} and contains results overlapping with some of the
work of \cite{PaulJ2007} and this paper.

The rest of the paper is organized as follows. In Section \ref{sec:model}, we
describe the model for the eigenvectors and analyze the risk of the standard
PCA estimator. In Section \ref{sec:lower_bound}, we present the lower bounds on
the minimax risk of any eigenvector estimator. In Section
\ref{sec:diagonal_thresh}, we derive a lower bound on the risk of the diagonal
thresholding estimator proposed by \cite{JohnstoneL2009}. In Section
\ref{sec:two_stage}, we propose a new estimator named ASPCA (augmented sparse
PCA) that is a refinement of the diagonal thresholding estimator.
In Section \ref{sec:summary}, we discuss the question of attainment of the risk
bounds. Proofs of the results are given in Section \ref{sec:proofs} in the
Appendix.

\begin{table}[!t]
\begin{center}
  \begin{tabular}{| l | c | c | c | }
    \hline
    {\bf Estimator} & {\bf dense} & {\bf sparse} & {\bf ultra-sparse} \\ \hline
    {\bf Lower bound} & $O(N/n)$  & $O(n^{-(1-q/2)})$ & $O((\log N/n)^{1-q/2})$\\ \hline
    PCA & rate optimal & inconsistent & inconsistent \\ \hline
    D.T. & inconsistent & not rate optimal & not rate optimal \\ \hline
    ASPCA & inconsistent  & ?  & rate optimal \\
    \hline
  \end{tabular}
  \caption{Comparison of Lower Bounds on eigenvector estimation and Worst Case Rates of various procedures.}
\end{center}
        \label{T:comparison}
\end{table}

\section{Problem setup}\label{sec:model}

First we introduce certain notations. Throughout, $\mathbb{S}^{N-1}$ denotes
the unit sphere in $\mathbb{R}^N$ centered at the origin, $\lfloor x \rfloor$
denotes the largest integer less than or equal to $x \in \mathbb{R}$.

Let  $\{X_i : i=1,\ldots,n\}$ be a triangular array, where for each
$n$, the $N \times 1$ random vectors $X_i := X_{i}^n, i=1,\ldots,n$ are
independent and identically distributed on a common probability space.
Throughout we assume that $X_i$'s are i.i.d. as
$N(\bs{0},\Sigma)$, where the population matrix $\Sigma$ is a finite rank
perturbation of (a multiple of) the identity. In other words,
\begin{equation}\label{eq:Sigma_basic}
\Sigma = \sum_{\nu=1}^M \lambda_\nu \theta_\nu \theta_\nu^T +
\sigma^2 I,
\end{equation}
where $\lambda_1 > \lambda_2 > \ldots > \lambda_M > 0$, and the
vectors $\theta_1,\ldots,\theta_M$ are orthonormal, which implies
(*).
$\theta_\nu$ is the eigenvector of $\Sigma$ corresponding to the
$\nu$-th largest eigenvalue, namely, $\lambda_\nu + \sigma^2$. The
term ``finite rank'' means that $M$ remains fixed even as $n \to
\infty$. The asymptotic setting involves letting both $n$ and $N$
grow to infinity simultaneously. For simplicity, we assume that the
$\lambda_\nu$'s are fixed while the parameter space for the
$\theta_\nu$'s varies with $N$.

The observations can be described in terms of the model
\begin{equation}\label{eq:basic}
X_{ik} = \sum_{\nu=1}^M \sqrt{\lambda_\nu} v_{\nu i} \theta_{\nu k} + \sigma
Z_{ik}, \quad i=1,\ldots,n, \quad k=1,\ldots,N.
\end{equation}
Here, for each $n$, $v_{\nu i}$, $Z_{ik}$ are i.i.d. $N(0,1)$. Since the
eigenvectors of $\Sigma$ are invariant to a scale change in the original
observations, it is assumed that $\sigma = 1$. Hence,
$\lambda_1,\ldots,\lambda_M$ in the asymptotic results should be replaced by
$\lambda_1/\sigma^2, \ldots, \lambda_M/\sigma^2$ when (\ref{eq:Sigma_basic})
holds with an arbitrary $\sigma > 0$. Since the main focus of this paper is
estimation of eigenvectors, without loss of generality we consider the
uncentered sample covariance matrix  $\mathbf{S} :=
\frac{1}{n}\mathbf{X}\mathbf{X}^T$, where $\mathbf{X} = [X_1:\ldots:X_n]$.

The following condition, termed \textit{Basic Assumption}, will be
used throughout the asymptotic analysis, and will be referred to as
{\bf BA}.
\begin{itemize}
\item[{\bf BA}]
(\ref{eq:basic}) holds with $\sigma = 1$; $N = N(n) \to \infty$ as $n \to
\infty$; $\lambda_1 > \ldots > \lambda_M
> 0$ are fixed (do not vary with $N$), where $M$ is unknown but fixed.
\end{itemize}

\subsection{Eigenvector estimation with squared error loss}\label{subsec:highdpca-lossfn}

Given data $\{X_i\}_{i=1}^n$,  the goal is to estimate $M$ and the eigenvectors
$\theta_1,\ldots, \theta_M$. For simplicity, to derive the lower
bounds, we first assume that $M$ is known. In Section
\ref{subsec:aspca-M_hat}
we derive an estimator of $M$, which can be shown to be consistent under the assumed sparsity
conditions. To assess the performance of any estimator, a minimax risk
analysis approach is proposed.
The first task is to specify a loss function $L(\widehat\theta_\nu,\theta_\nu)$
between the estimated and true eigenvector. Since the model is invariant to sign changes of each $\theta_\nu$, we
consider the following loss function, also invariant to sign changes.
\begin{equation}\label{eq:loss_fn}
L(\mathbf{a},\mathbf{b}) := 2(1 - |\langle \mathbf{a},\mathbf{b}\rangle|) =
\parallel \mathbf{a} - sign(\langle \mathbf{a}, \mathbf{b} \rangle) \mathbf{b}
\parallel^2,
\end{equation}
where $\mathbf{a}$ and $\mathbf{b}$ are $N \times 1$ vectors with unit $l_2$
norm. An estimator $\widehat\theta_\nu$ is called consistent with respect to
\(L\), if $L(\widehat \theta_\nu, \theta_\nu) \to 0$ in probability as $n \to
\infty$.
%


\subsection{Rate of convergence for ordinary PCA}

We first consider the asymptotic risk of the
leading eigenvectors of the sample covariance matrix (henceforth
referred to as the standard PCA estimators) when the ratio $N/n$ is
small. Specifically, it is assumed that $N/n \to 0$ as $n \to
\infty$.

For future use, we define
\begin{equation}\label{eq:h_lambda}
h(\lambda) := \frac{\lambda^2}{1+\lambda}  \qquad \lambda > 0,
\end{equation}
and
\begin{equation}\label{eq:g_lambda_tau}
g(\lambda,\tau) = \frac{(\lambda-\tau)^2}{(1+\lambda)(1+\tau)},
\qquad \lambda, \tau
> 0.
\end{equation}
In \cite{JohnstoneL2009} (Theorem 1) it was shown that under a single spike model, as $N/n \to
0$, the standard PCA estimator of the leading eigenvector is consistent. The following result, proven in the Appendix, is a
refinement of that, as it also provides the leading error term.

\begin{theorem}\label{thm:OPCA_risk_bound}
Let $\widehat\theta_{\nu,PCA}$ be the eigenvector corresponding to the $\nu$-th
largest eigenvalue of $\mathbf{S}$. Assume that {\bf BA} holds and $N,n \to
\infty$ such that $N/n \to 0$, and moreover, $\log n = o(N)$. Then, for each
$\nu=1,\ldots,M$,
\begin{equation}\label{eq:OPCA_risk_bound}
\sup_{\theta_\nu \in \mathbb{S}^{N-1}} \mathbb{E}L(\widehat
\theta_{\nu,PCA}, \theta_\nu) =
    \left[\frac{N-M}{nh(\lambda_\nu)} + \frac{1}{n}\sum_{\mu \neq \nu}
\frac{1}{g(\lambda_\mu,\lambda_\nu)}\right](1+o(1)).
\end{equation}
\end{theorem}

\begin{remark}
Observe that Theorem \ref{thm:OPCA_risk_bound} does not assume any special
structure (e.g., sparsity) for the eigenvectors. The first term on the RHS of
(\ref{eq:OPCA_risk_bound}) is a nonparametric component which arises from the
interaction of the noise terms with the different coordinates, while the second
term is a parametric component which results from the interaction with the
remaining $M-1$ eigenvectors corresponding to different eigenvalues. The second
term shows that the closer the successive eigenvalues are, the larger is the
estimation error. The upshot of (\ref{eq:OPCA_risk_bound}) is that standard PCA
provides a consistent estimator of the leading eigenvectors of the population
covariance matrix when the dimension-to-sample-size ratio ($N/n$) is
asymptotically negligible.
\end{remark}

\subsection{$l_q$ constraint on eigenvectors}\label{subsec:l_q_constraint}

As shown by various authors \citep{Nadler2008,Onatski2006, Paul2007}, when $N/n \to c \in (0,\infty]$, standard PCA provides inconsistent
estimators for the population eigenvectors. In this subsection we consider the
following model for approximate sparsity of the eigenvectors.
For each $\nu=1,\ldots,M$, we assume that $\theta_\nu$ belongs to an $l_q$ ball
with radius $C$, for some $q \in (0,2)$. Specifically, we assume that
$\theta_\nu \in \Theta_q(C)$, where
\begin{equation}\label{eq:Theta_q_C}
\Theta_q(C) := \{ \bs{a} \in \mathbb{S}^{N-1} : \sum_{k=1}^N |a_k|^q \leq
C^q\}.
\end{equation}
Note that our condition of sparsity is slightly different from that of
\cite{JohnstoneL2009}.


Note that since $0< q < 2$, for $\Theta_q(C)$ to be
nonempty, one needs $C \geq 1$.
Further, if $C^q \geq N^{1-q/2}$, then the space $\Theta_q(C)$ is all of
$\mathbb{S}^{N-1}$ because in this case, the least sparse vector
$\frac1{\sqrt{N}}(1,1,\ldots,1)$ is in the parameter space.


The parameter space for $\boldsymbol{\theta} :=
[\theta_1:\ldots:\theta_M]$ is denoted by
\begin{equation}\label{eq:Theta_q_M}
\Theta_q^M(C_1,\ldots,C_M) := \{\boldsymbol{\theta} \in
\prod_{\nu=1}^M \Theta_q(C_\nu) ~:~ \langle \theta_\nu,
\theta_{\nu'} \rangle = 0, ~~\mbox{for}~~ \nu \neq \nu'\},
\end{equation}
where $\Theta_q(C)$ is defined through (\ref{eq:Theta_q_C}), and
$C_\nu \geq 1$ for all $\nu=1,\ldots,M$.

\begin{remark}\label{rem:covariance_sparsity}
While our focus is on eigenvector sparsity, condition (\ref{eq:Theta_q_M}) also
implies sparsity of the covariance matrix itself. In particular, for $q \in
(0,1)$, a spiked covariance matrix satisfying (\ref{eq:Theta_q_M}) also belongs
to the class of sparse covariance matrices analyzed by \cite{BickelL2008b},
\cite{CaiL2011} and \cite{CaiZ2011}. Indeed, \cite{CaiZ2011} obtained the
minimax rate of convergence for covariance matrix estimators under the spectral
norm when the rows of the population matrix satisfy a weak-$l_q$ constraint.
However, as we will show below, the minimax rate for estimation of the leading
eigenvectors is faster than that for covariance estimation.
%
%
\end{remark}

\section{Lower bounds on the minimax risk}\label{sec:lower_bound}

We now derive lower bounds on the minimax risk of estimating $\theta_\nu$ under
the loss function (\ref{eq:loss_fn}). To aid in describing and interpreting the
lower bounds, we define the following two auxiliary parameters.
The first is an \textit{effective noise level per coordinate}
\begin{equation}
  \label{eq:tau-nu2}
  \tau_\nu^2 = 1/(n h(\lambda_\nu))
\end{equation}
and the second is an \textit{effective dimension}
\begin{equation}
  \label{eq:m-nu}
  m_\nu := A_q (\bar C_\nu/\tau_\nu)^q
\end{equation}
where $a_q := (2/9)^{1-q/2}$, $c_1 := \log(9/8)$ and $A_q := 1/(a_q c_1^{q/2})$
and $\bar C_\nu^q := C_\nu^q - 1$.

The phrase \textit{effective noise level per coordinate} is motivated by the
risk bound in Theorem \ref{thm:OPCA_risk_bound}, since dividing both sides of
(\ref{eq:OPCA_risk_bound}) by $N$, the expected ``per coordinate'' risk (or
variance) of the PCA estimator is asymptotically $\tau_\nu^2$.
Next, following \cite{Nadler2009}, let us provide a different interpretation of
$\tau_\nu$. Consider a sparse $\theta_\nu$ and an oracle that, regardless of
the observed data, selects a set $J_\tau$ of all coordinates of $\theta_\nu$
that are larger than $\tau$ in absolute value, and then performs PCA on the
sample covariance restricted to these coordinates. Since $\theta_\nu \in
\Theta_q(C_\nu)$, the maximal squared-bias is
\begin{eqnarray*}
\sup_{\theta_\nu \in \Theta_q(C_\nu)} \sum_{k \not\in J_\tau} |\theta_{\nu
k}|^2 &\asymp& \sup\{\sum_{k=1}^N x_k^{2/q} : \sum_{k=1}^N x_k \leq
C_\nu^q, \max_{k} x_k < \tau^q, \min_{k} x_k \geq 0  \} \\
&\asymp& C_\nu^q \tau^{2-q}
\end{eqnarray*}
which follows by the correspondence $x_k = |\theta_{\nu k}|^q$, and  the
convexity of the function $\sum_{k=1}^N x_k^{2/q}$. On the other hand, by
Theorem \ref{thm:OPCA_risk_bound}, the maximal variance term of this oracle
estimator is of the order $k_\tau/(nh(\lambda_\nu))$ where $k_\tau$ is the
maximal number of coordinates of $\theta_\nu$ exceeding $\tau$. Again,
$\theta_\nu \in \Theta_q(C_\nu)$ implies that $k_\tau \asymp C_\nu^q
\tau^{-q}$. Thus, to balance the bias and variance terms, we need $\tau \asymp
1/\sqrt{nh(\lambda_\nu)} = \tau_\nu$. This heuristic analysis shows that
$\tau_\nu$ can be viewed as an \textit{oracle threshold} for the coordinate
selection scheme, i.e., the best possible estimator of $\theta_\nu$ based on
individual coordinate selection can expect to recover only those coordinates
that are above the threshold $\tau_\nu$.

To understand why $m_\nu$ is an \textit{effective dimension}, consider the
least sparse vector $\theta_\nu\in \Theta_q(C_\nu)$.
This vector should have as many nonzero coordinates of equal size as possible.
If $C_\nu^q > N^{1-q/2}$ then the vector with coordinates $\pm N^{-1/2}$ does
the job. Otherwise, we set the first coordinate of the vector to be
$\sqrt{1-r^2}$ for some $r \in (0,1)$ and choose all the nonzero coordinates to
be of magnitude $\tau_\nu$. Clearly, we must have $r^2 = m \tau_\nu^2$, where
$m+1$ is the maximal number of nonzero coordinates, while the $l_q$ constraint
implies that $(1-r^2)^{q/2} + m \tau_\nu^q \leq C_\nu^q$. The last inequality
shows that the maximal $m$ is just a constant multiple of $m_\nu$. This
construction also constitutes the key idea in the proof of Theorems
\ref{th:three-way-lower} and \ref{th:sparse-lower}. Finally, we set
\begin{equation}
  \label{eq:Nprime}
  N' = c_1(N-M),
\end{equation}
where the origin of $c_1 = \log(9/8)$ will be explained in the
proof.


\begin{theorem}
\label{th:three-way-lower}
Assume that {\bf BA} holds, $0 < q < 2$, and $n, N \to \infty$.
Then, there exists a constant $B_1 > 0$ such that for $n$
sufficiently large,
  \begin{equation}
    \label{eq:mmxbd}
    R_\nu^* := \inf_{\widehat \theta_\nu} \sup_{\Theta_q(\mathbf{C}) } \mathbb{E}
    L(\widehat \theta_\nu, \theta_\nu) \geq B_1 \delta_n,
  \end{equation}
where $\delta_n$ is given by
\begin{equation*}
  \delta_n =
  \left\{
  \begin{array}{cll}
     \tau_\nu^2 N'    & \quad \text{if} \quad \tau_\nu^2  N' < 1 ~\text{and}~ N' < m_\nu   & ~~~[\mbox{dense setting}]\\
     \tau_\nu^2 m_\nu & \quad \text{if} \quad \tau_\nu^2  m_\nu < 1 ~\text{and}~m_\nu < N' & ~~~[\mbox{sparse setting}]\\
     ~~1  & \quad \text{if} \quad \tau_\nu^2 \cdot \min\{N',m_\nu\}  > 1 & ~~~[\mbox{weak signal}].
  \end{array}
  \right.
\end{equation*}
\end{theorem}

We may think of $m_n := \min\{N', m_\nu\}$ as the effective dimension of the
least favorable configuration. In the \textit{sparse} setting, $m_n = A_q
\bar{C}_\nu^q [n h(\lambda_\nu)]^{q/2} < c_1 N$ (i.e., $\bar{C}_\nu^q n^{q/2} <
c' N$ for some $c' > 0$), and the lower bound is of the order
\begin{equation}\label{eq:sparse_rate}
  \delta_n = c_1 A_q C_\nu^q \tau_\nu^{2-q} = \frac{ c_1 A_q C_\nu^q}{[n h(\lambda_\nu)]^{1 - q/2}}
            \asymp \frac{C_\nu^q}{n^{1-q/2}}~.
\end{equation}
On the other hand, in the \textit{dense} setting, $m_n = c_1(N-M)$. If $N/n \to
c$ for some $c > 0$, then $\delta_n =c_1(N-M)/(n h(\lambda_\nu)) \asymp 1$,
and so any estimator of the eigenvector $\theta_\nu$ is
inconsistent. If $N/n \to 0$ then the lower bound is
\begin{equation}\label{eq:dense_rate}
  \delta_n = \frac{c_1(N-M)}{n h(\lambda_\nu)} \asymp
  \frac{N}{n}~.
\end{equation}
%
Eq. (\ref{eq:dense_rate}) and Theorem \ref{thm:OPCA_risk_bound}
imply that in the dense setting with $N/n\to 0$, the standard PCA estimator
$\widehat\theta_{\nu,PCA}$ attains the optimal rate of convergence.

A sharper lower bound is possible in what we call an
\textit{ultra-sparse} setting which happens if $\bar{C}_\nu^q
n^{q/2} = O(N^{1-\alpha})$ for some $\alpha \in (0,1)$. In
this case the dimension $N$ is much larger than the quantity
$\bar{C}_\nu^q n^{q/2}$ measuring the effective dimension. Hence,
we define a modified effective noise level
per-coordinate
\begin{displaymath}
  \bar \tau_\nu^2 = \frac{\alpha}{9} \frac{\log N}{n h(\lambda_\nu)},
\end{displaymath}
and a modified effective dimension
\begin{displaymath}
  \bar m_\nu = a_q^{-1} (\bar C_\nu/ \bar \tau_\nu)^q.
\end{displaymath}

\begin{theorem}\label{th:sparse-lower}
  Assume that {\bf BA} holds, $0 < q < 2$, and $n,N \to \infty$ such that $\bar m_\nu
  = O(N^{1 - \alpha})$ for some $\alpha \in (0,1)$. Then, assuming that
  $\bar m_\nu \bar \tau_\nu^2 \leq 1$
  for $n$ sufficiently large, the
  minimax bound \eqref{eq:mmxbd} holds with
  \begin{equation}
            \label{eq:ultrasparse_rate}
    \delta_n = \bar m_\nu \bar \tau_\nu^2
          = a_q^{-1} C_\nu^q \Big( \frac{ \log N}{n h(\lambda_\nu)}
          \Big)^{1 - q/2}. \quad
          \mbox{[ultra-sparse setting]}
  \end{equation}
\end{theorem}
Note that in the ultra-sparse setting $\delta_n$ is larger by a factor of $(\log N)^{1-q/2}$
compared to the sparse setting, Eq. (\ref{eq:sparse_rate}).

\section{Risk of the diagonal thresholding
estimator}\label{sec:diagonal_thresh}

In this section, we analyze the convergence rate of the SPCA scheme (henceforth
referred to as the diagonal thresholding or D.T. scheme) proposed by
\cite{JohnstoneL2009}. In this section and in Section \ref{sec:two_stage}, we
assume for simplicity that $N \geq n$.
Let the sample variance of the $k$-th
coordinate (i.e., the $k$-th diagonal entry of $\mathbf{S}$) be
denoted by $\mathbf{S}_{kk}$. Then the D.T. scheme consists of the
following steps.
\begin{enumerate}
\item
Define $I=I(\gamma_n)$ to be the set of indices $k \in
\{1,\ldots,N\}$ such that $\mathbf{S}_{kk} > \gamma_n$ for some
threshold $\gamma_n > 0$.
\item
Let $\mathbf{S}_{II}$ be the submatrix of $\mathbf{S}$ corresponding
to the coordinates $I$. Perform an eigen-analysis of
$\mathbf{S}_{II}$. Denote the eigenvectors by
$\mathbf{f}_1,\ldots,\mathbf{f}_{\min\{n,|I|\}}$.
\item
For $\nu=1,\ldots,M$, estimate $\theta_\nu$ by
the $N\times 1$ vector
$\widetilde
{\mathbf{f}}_\nu$, obtained from $\mathbf{f}_\nu$ by augmenting zeros to
all the coordinates in $I^c := \{1,\ldots,N\} \setminus I$.
\end{enumerate}

Assuming that $\theta_\nu\in \Theta_q(C_\nu)$, \cite{JohnstoneL2009}
showed that the D.T. scheme with a threshold of the form
$\gamma_n = 1 + \gamma \sqrt{\log N/n}$ for some $\gamma > 0$
leads to a consistent estimator of $\theta_\nu$.
The risk of this estimator, however, was not analyzed in \cite{JohnstoneL2009}.
As we prove below, the risk of the D.T. estimator is not rate optimal.
This can be anticipated from the lower bound on the minimax risk (Theorems
\ref{th:three-way-lower} and \ref{th:sparse-lower}) which indicate
that to attain the optimal risk, a coordinate selection
scheme must select all coordinates of $\theta_\nu$ of size at least
$c \sqrt{\log N/n}$. With a threshold of the form $\gamma_n$ above, however,
only coordinates of size $(\log N/n)^{1/4}$ are selected.
As shown in the following theorem, even for the case of a single signal ($M=1$)
this leads to a much larger lower bound.

\begin{theorem}\label{thm:diagonal_thresholding_risk}
Suppose that {\bf BA} holds with $M =1$. Let $C > 0$, $0 < q < 2$, and $n,N \to
\infty$ be such that $C^q n^{q/4} = o(\max\{\sqrt{n},N\})$. Then the Diagonal
Thresholding estimator $\widehat \theta_{1,DT}$ proposed by
\cite{JohnstoneL2009} satisfies, for any $q \in (0,2)$,
\begin{equation}
    \label{eq:DT_risk}
\sup_{\theta_1 \in \Theta_q(C)}
\mathbb{E}L(\widehat\theta_{1,DT},\theta_1) \geq K_q \bar{C}^q
n^{-\frac{1}{2}(1-q/2)}
\end{equation}
for a constant $K_q > 0$, where $\bar{C}^q = C^q -1$.
\end{theorem}
Comparing (\ref{eq:DT_risk}) with the lower bound (\ref{eq:sparse_rate}), shows
the large gap between the two rates, $n^{-1/2(1-q/2)}$ vs. $n^{-(1-q/2)}$. The
reason for this difference is that the D.T. scheme uses only the diagonal
of the sample covariance matrix \(\bf S\), ignoring the information in its off-diagonal
entries. In the next section we propose a refinement of the D.T. scheme,
denoted ASPCA,
that constructs an improved eigenvector estimate using all entries of \(\bf S\).

\section{A two stage coordinate selection scheme}\label{sec:two_stage}


As discussed above, the DT scheme can reliably detect only those eigenvector
coordinates $|\theta_{\nu,k}| = O((\log N/n)^{1/4})$, whereas to reach the lower
bound one needs to detect those coordinates of size $|\theta_{\nu,k}| = O((\log N/n)^{1/2})$.

To motivate an improved coordinate selection scheme, consider
a partition of the
$N$ coordinates into two sets $A$ and $B$, where the
former contains all those $k$ such that $|\theta_{1k}|$ is ``large''
(selected by the D.T. scheme), and the latter contains the remaining
smaller coordinates. Partition the matrix $\Sigma$ as
\begin{equation*}
\Sigma = \begin{bmatrix} \Sigma_{AA} & \Sigma_{AB} \cr
                         \Sigma_{BA} & \Sigma_{BB} \cr
         \end{bmatrix}.
\end{equation*}
Observe that, $\Sigma_{BA} = \lambda_1 \theta_{1,B} \theta_{1,A}^T$. Let
$\widetilde \theta_1$ be a ``preliminary'' estimator of $\theta_1$ such that
$\lim_{n \to \infty} \mathbb{P}(\langle \widetilde \theta_{1,A},
\theta_{1,A}\rangle \geq \delta_0) = 1$ for some $\delta_0 > 0$ (e.g.,
$\widetilde \theta_1$ could be the D.T. estimator).
Then we have the relationship,
\begin{equation*}
\Sigma_{BA} \widetilde \theta_{1,A} = \langle \widetilde \theta_{1,A},
\theta_{1,A}\rangle \lambda_1 \theta_{1,B} \approx c(\delta_0) \lambda_1
\theta_{1,B}
\end{equation*}
for some $c(\delta_0)$ bounded below by $\delta_0/2$, say. Thus, one possible
strategy is to additionally select all those coordinates of $\Sigma_{BA}
\widetilde \theta_{1,A}$ that are larger (in absolute value) than some constant
multiple of $\sqrt{\log N}/\sqrt{nh(\lambda_1)}$. In practice we do not know
$\Sigma_{BA}$ or $\lambda_1$ but we can use $\mathbf{S}_{BA}$ as a surrogate
for the former and the largest eigenvalue of $\mathbf{S}_{AA}$ to obtain an
estimate for the latter. A technical challenge is to show, that with
probability tending to 1, such a scheme indeed recovers all coordinates $k$
with $|\theta_{1 k}| > c_1 \sqrt{\log N}/\sqrt{nh(\lambda_1)}$, while
discarding all coordinates $k$ with $|\theta_{1k}| < c_2  \sqrt{\log
N}/\sqrt{nh(\lambda_1)}$ for some constants $c_1 > c_2
> 0$. Figure 1 provides a pictorial description of
the D.T. and ASPCA coordinate coordinate selection schemes.

\begin{figure}[h]\label{fig:thresholding_schemes}
\begin{center}
\includegraphics[width=5in, height=4in, angle=0]{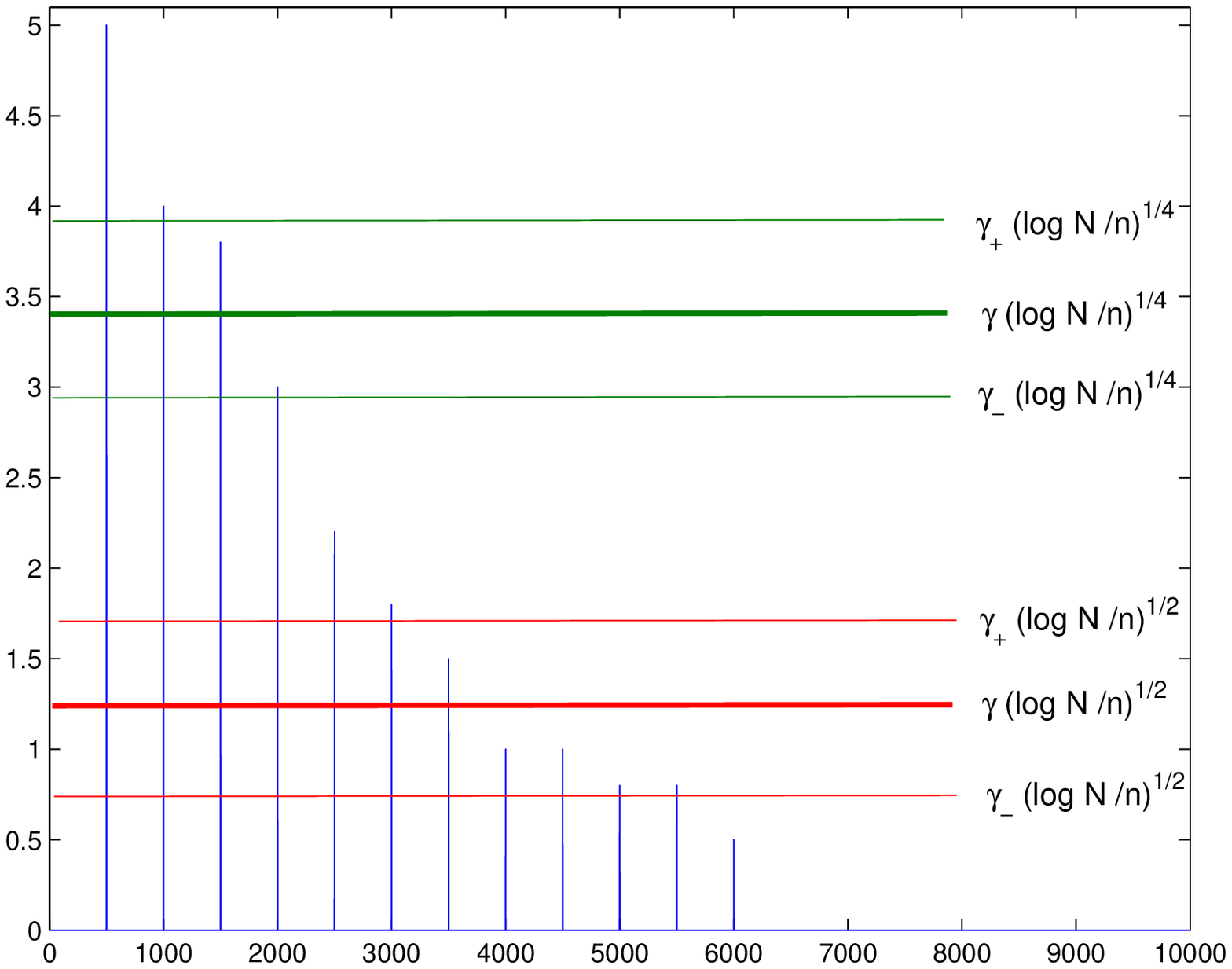}
\caption{Schematic diagram of the D.T. and ASPCA thresholding schemes under the
single component setting. The vertical lines depict the absolute values of the
coordinates of the first eigenvector. The threshold for the D.T. scheme is
$\gamma (\log N/n)^{1/4}$ while the thresholds for the ASPCA scheme is
$\gamma(\log N/n)^{1/2}$. The schemes select the coordinates above the upper
limits (indicated by the multiplier $\gamma_+$) and discard the coordinates
below the lower limits (indicated by multiplier $\gamma_-$) with high
probability. Here, $\gamma_+ > \gamma > \gamma_- > 0$ are generic constants.}
\end{center}
\end{figure}

\subsection{ASPCA scheme}

Based on the ideas described above, we now present the ASPCA algorithm.
It first makes two stages of coordinate selection, whereas
the final stage consists of an eigen-analysis of the
submatrix of $\mathbf{S}$ corresponding to the selected coordinates.
The algorithm is described below.

For any $\gamma > 0$ define
\begin{equation}\label{eq:I_gamma_def}
I(\gamma) = \{k : \mathbf{S}_{kk} > 1 +\gamma\}.
\end{equation}
Let $\gamma_i > 0$ for $i=1,2$ and $\kappa > 0$  be constants to be
specified later.
\begin{itemize}
\item[{\bf Stage 1}]
\item[{$1^o$}] Let
$I = I(\gamma_{1,n})$ where $\gamma_{1,n} = \gamma_1 \sqrt{\log
N/n}$.
\item[{$2^o$}]
Denote the eigenvalues and eigenvectors  of $\mathbf{S}_{II}$ by
$\widehat\ell_1 > \ldots > \widehat\ell_{m_1}$  and
$\mathbf{f}_1,\ldots,\mathbf{f}_{m_1}$ respectively, where $m_1 =
\min\{n, | I|\}$,
\item[{$3^o$}]
Estimate $M$ by $\widehat M$ defined in Section
\ref{subsec:aspca-M_hat}.

\item[{\bf Stage 2}]
\item[{$4^o$}] Let
$\mathbf{E} = [\widehat\ell_1^{-1/2}  \mathbf{f}_1 \cdots
\widehat\ell_{\widehat
  M}^{-1/2} \mathbf{f}_{\widehat M} ]$ and
 $ \mathbf{Q} = \mathbf{S}_{I^c I} \mathbf{E}.$
\item[{$5^o$}]
Let $J = \{ k \not\in I ~:~ (\mathbf{Q}\mathbf{Q}^T)_{kk} >
\gamma_{2,n}^2\}$ for some $\gamma_{2,n} > 0$. Define  $K = I \cup
J$.
\item[{\bf Stage 3}]
\item[{$6^o$}]
For $\nu=1,\ldots, \widehat M$, denote by $\widehat \theta_\nu$ the
$\nu$-th eigenvector of $\mathbf{S}_{KK}$, augmented with zeros in
the coordinates $K^c$.
\end{itemize}

\begin{remark}\label{rem:gamma_choice}
The ASPCA scheme is specified up to the choice of parameters
$\gamma_1,\gamma_{2,n}$ and $\kappa$, that determine its rate of convergence.
It can be shown that choosing $\gamma_1 = 4$, $\kappa = \sqrt{2+\epsilon}$ for
some $\epsilon>0$, and $\gamma_{2,n}$ given by
\begin{equation}\label{eq:gamma_2_n_choice}
\gamma_{2,n} = \gamma_2 \left(\sqrt{\frac{\log N}{n} } +
\frac{1}{\kappa} \sqrt{\frac{\widehat M}{n}}\right)
\end{equation}
with $\gamma_2 = \kappa \sqrt{3/2}$ results in an asymptotically optimal rate.
Again, we note that for finite $N$, $n$, the actual performance in terms of the
risk of the resulting eigenvector estimate may have a strong dependence on the
threshold. In practice, a delicate choice of thresholds can be highly
beneficial. This issue, as well as the analysis of the risk of the ASPCA
estimator, are beyond the scope of this paper and will be studied in a separate
publication.
\end{remark}

\subsection{Estimation of $M$}\label{subsec:aspca-M_hat}


Estimation of the dimension of the signal subspace is a classical
problem. If the signal eigenvalues are strong enough (i.e.,  $\lambda_\nu
> c \sqrt{N/n}$ for all $\nu=1,\ldots,M$, for some $c > 1$ independent of $N,n$), then
nonparametric methods that do not assume eigenvector sparsity can
asymptotically estimate the correct $M$ (see, e.g.
\cite{Kritchman2008}). When the eigenvectors are sparse, we can
detect much weaker signals, as we describe below.

We estimate $M$ by thresholding the eigenvalues of the submatrix
$\mathbf{S}_{\bar I \bar I}$ where $\bar I := I(\bar\gamma\sqrt{\log
N/n})$ for some $\bar\gamma > 0$. Let $\bar m = \min \{ n, |\bar I|
\}$ and $\bar \ell_1 > \ldots > \bar \ell_{\bar{m}}$ be the nonzero
eigenvalues of $\mathbf{S}_{\bar I \bar I}$. Let $\alpha_n
> 0$ be a user-defined threshold. Then, define $\widehat M$ by
\begin{equation}\label{eq:M_hat_def}
\widehat M := \max \{1\leq k \leq \bar m : \bar
  \ell_k  > 1+\alpha_n \}.
\end{equation}
It can be shown that under appropriate sparsity conditions, with a suitable
choice of threshold $\alpha_n$, $\widehat M$ is a consistent estimator of $M$.


\section{Summary and Discussion}\label{sec:summary}

In this paper we derived lower bounds on eigenvector estimates under three
different sparsity regimes, denoted dense, sparse, and ultra-sparse. In the
\textit{dense} setting, Theorems \ref{thm:OPCA_risk_bound} and
\ref{th:three-way-lower} show that when $N/n \to 0$, the standard PCA estimator
attains the optimal rate of convergence. In the \textit{ultra-sparse} setting,
Theorem 3.1 of \cite{Ma2011} shows that the maximal risk of the ITSPCA
estimator proposed by him attains the same asymptotic rate as the corresponding
lower bound of Theorem \ref{th:sparse-lower}. This implies that in the
ultra-sparse setting, the lower bound on the minimax rate is indeed sharp. In
a separate paper, we prove that  in the ultra-sparse regime, the ASPCA
algorithm also attains the minimax rate.

Finally, our analysis leaves some open questions in the intermediate sparse
regime. According to Theorem \ref{th:three-way-lower}, the lower bound in this
regime is smaller by a factor of \((\log N)^{1-q/2}\), as compared to the
ultra-sparse setting. Therefore, whether there exists an estimator (and in
particular, one with low complexity), that attains the current lower bound, or
whether this lower bound can be improved is an open question for future
research.

\appendix

\setcounter{section}{0}

\setcounter{equation}{0}

\setcounter{lemma}{0}

\setcounter{lemma}{0}

\setcounter{proposition}{0}

\setcounter{remark}{0}

\renewcommand{\theequation}{A.\arabic{equation}}
\renewcommand{\thetable}{A.\arabic{table}}
\renewcommand{\thefigure}{A.\arabic{figure}}
\renewcommand{\thelemma}{A.\arabic{lemma}}
\renewcommand{\theproposition}{A.\arabic{proposition}}
\renewcommand{\theremark}{A.\arabic{remark}}

\section{Proofs}\label{sec:proofs}

\subsection{Asymptotic risk of the standard PCA estimator}\label{sec:OPCA_risk_bound}

To prove Theorem \ref{thm:OPCA_risk_bound}, on the risk of the PCA\ estimator, we use the following lemmas.


\subsubsection*{Deviation of extreme eigenvalues of Wishart matrices}

In our analysis, we shall need a probabilistic bound for deviations of
$\parallel \frac{1}{n} \mathbf{Z}\mathbf{Z}^T - I\parallel$. This is given
in the following lemma, proven in Section \ref{sec:auxiliary_results}.
\begin{lemma}\label{lemma:eigen_deviation_bound}
Let $t_n = 8(N_n/n)\sqrt{\log N_n / N_n}$ where $N_n = \max\{n, N\}$. Let
$\mathbf{Z}$ be an $N \times  n$ matrix with i.i.d. $N(0,1)$ entries. Then for
any $c
> 0$, there exists $n_c\geq 1$ such that for all $n \geq n_c$,
\begin{equation}
\mathbb{P}\left( \parallel \frac{1}{n} \mathbf{Z}\mathbf{Z}^T - I_N
\parallel  > \frac{N}{n} + 2\sqrt{\frac{N}{n}} + ct_n \right) \leq 2
N_n^{-c^2}.
\end{equation}
\end{lemma}

\subsubsection*{Deviation of quadratic
forms}\label{subsec:deviation_quad_form}

The following lemma is due to \cite{Johnstone2001b}.
\begin{lemma} \label{lemma:chi_square_large_dev}
Let $\chi_{n}^2$ denote a Chi-square random variable with $n$ degrees of
freedom. Then,
\begin{eqnarray}
 \mathbb{P}(\chi_{n}^2 > n(1+\epsilon) ) &\leq& e^{-3n\epsilon^2/16}
\qquad (0 < \epsilon < \frac{1}{2}),\label{eq:large_dev_chisq_1}\\
\mathbb{P}(\chi_{n}^2 < n(1-\epsilon) ) &\leq& e^{-n\epsilon^2/4}
~~~~~~~~ (0 < \epsilon < 1),\label{eq:large_dev_chisq_2}\\
\mathbb{P}(\chi_{n}^2 > n(1+\epsilon) ) &\leq& \frac{\sqrt{2}}{\epsilon
  \sqrt{n}}e^{-n\epsilon^2/4}
~~(0 < \epsilon < 1/2, n \geq 16). \label{eq:large_dev_chisq_3}
\end{eqnarray}
\end{lemma}

\noindent The following lemma is from \cite{JohnstoneL2009}.
\begin{lemma} \label{lemma:quad_form_large_dev}
Let  $y_{1i},y_{2i},i=1,\ldots,n$ be two sequences of mutually
  independent, i.i.d. $N(0,1)$ random variables. Then for large $n$
  and any $b$ s.t. $0 < b \ll \sqrt{n}$,
\begin{equation}\label{eq:large_dev_covar}
\mathbb{P}\left(|\frac{1}{n} \sum_{i=1}^n y_{1i} y_{2i}|>\sqrt{b/n}\right) \leq
2\exp\left\{-\frac{3b}{2} + O(n^{-1}b^2)\right\}.
\end{equation}
\end{lemma}

\subsubsection*{Perturbation of
  eigen-structure}\label{subsec:highdpca-perturb_eigen}

The following lemma from \cite{Paul2005} is convenient for risk
analysis of estimators of eigenvectors. Several variants of this
lemma appear in the literature,
most based on the approach of \cite{Kato1980}.

\begin{lemma} \label{lemma:evec_perturb_bound}
Let $A$ and $B$ be two symmetric $m \times m$ matrices. Let the
eigenvalues of matrix $A$ be denoted by $\lambda_1(A) \geq \ldots
\geq \lambda_m(A)$. Set $\lambda_0(A) = \infty$ and
$\lambda_{m+1}(A) = - \infty$. For any $r \in\{1,\ldots,m\}$, if
$\lambda_r(A)$ is a unique eigenvalue of $A$, i.e., if
$\lambda_{r-1}(A) > \lambda_r(A)
> \lambda_{r+1}(A)$, then denoting by $\mathbf{p}_r$ the eigenvector
associated with the $r$-th eigenvalue,
\begin{equation}\label{eq:eig_perturb_first}
\mathbf{p}_r(A+B) - \mbox{sign}(\mathbf{p}_r(A+B)^T\mathbf{p}_r(A))
\mathbf{p}_r(A) = - H_r(A) B \mathbf{p}_r(A) + R_r
\end{equation}
where $H_r(A) := \sum_{s \neq r} \frac{1}{\lambda_s(A) -
\lambda_r(A)} P_{{\cal E}_s}(A)$ and $P_{{\cal E}_s}(A)$ denotes the
projection matrix onto the eigenspace ${\cal E}_s$ corresponding to
eigenvalue $\lambda_s(A)$ (possibly multi-dimensional). Define
$\Delta_r$ and  $\overline{\Delta}_r$ as
\begin{eqnarray}\label{eq:eigen_Delta_bar_r}
\Delta_r &:=& \frac{1}{2} [\parallel H_r(A) B \parallel + |\lambda_r(A+B)
- \lambda_r(A)| \parallel H_r(A) \parallel] \label{eq:eigen_Delta_r} \\
\overline{\Delta}_r &=& \frac{\parallel B\parallel} {\min_{1\leq j \neq r \leq
m} |\lambda_j(A) - \lambda_r(A)|}~.
\end{eqnarray}
Then, the residual term $R_r$ can be bounded by
\begin{eqnarray}\label{eq:eigenvec_error}
&& \hskip-0.2in \parallel R_r \parallel ~\leq~ \min\left\{ 10\overline{\Delta}_r^2, \right. \nonumber\\
&& \hskip-0.2in  \left.\parallel H_r(A) B \mathbf{p}_r(A)
\parallel \left[\frac{2\Delta_r(1+2\Delta_r)}{1-2\Delta_r(1+2\Delta_r)} +
\frac{\parallel H_r(A) B \mathbf{p}_r(A)\parallel}
{(1-2\Delta_r(1+2\Delta_r))^2}\right]\right\}
\end{eqnarray}
where the second bound holds only if $\Delta_r < (\sqrt{5}-1)/4$.
\end{lemma}

\begin{remark}\label{rem:perturbation_bound}
We can simplify the bound on the perturbation in (\ref{eq:eigenvec_error}) to
show that if $\overline{\Delta}_r \leq 1/4$, then
\begin{equation}\label{eq:eigenvec_error_leading}
\parallel R_r \parallel \leq C \parallel H_r(A) B \mathbf{p}_r(A) \parallel \overline{\Delta}_r
\end{equation}
where we can take $C=30$. To see this, note that $|\lambda_r(A+B) -
\lambda_r(A)| \leq
\parallel B \parallel$ and that $\parallel H_r(A) \parallel \leq [\min_{j\neq r}|\lambda_j(A)
-\lambda_r(A)|]^{-1}$, so that,
$$
\Delta_r \leq \parallel H_r(A) \parallel \parallel B \parallel \leq
\overline{\Delta}_r.
$$
Now, defining $\delta := 2\overline{\Delta}_r
(1+2\overline{\Delta}_r)$ and $\beta := \parallel H_r(A) B
\mathbf{p}_r(A) \parallel$, we have $10 \overline{\Delta}_r^2 \leq
(5/2) \delta^2$, and the bound (\ref{eq:eigenvec_error}) may be
expressed as
$$
\parallel R_r \parallel \leq  \frac{\beta \delta}{1-\delta} \min \left\{ \frac{5}{2} \frac{\delta(1-\delta)}{\beta},
1+\frac{\beta}{\delta(1-\delta)} \right\}~.
$$
For $x > 0$, the function $x \mapsto \min\{5x/2, 1+1/x\} \leq 5/2$.
Further, if $\overline{\Delta}_r < 1/4$, then $\delta <
3\overline{\Delta}_r < 3/4$ and so we conclude that
$$
\parallel R_r \parallel  \leq 10 \beta \delta \leq 30 \beta \overline{\Delta}_r.
$$
\end{remark}

For notational simplicity, throughout this subsection, we write $\widehat
\theta_\nu$ to mean $\widehat\theta_{\nu,PCA}$. Recall that the loss function
$L(\widehat\theta_\nu,\theta_\nu) = \parallel \widehat\theta_\nu -
\mbox{sign}\langle \widehat\theta_\nu,\theta_\nu\rangle \theta_\nu
\parallel^2$.
Invoking Lemma \ref{lemma:evec_perturb_bound} with $A = \Sigma$ and
$B=\mathbf{S}-\Sigma$ we get
\begin{equation}\label{eq:theta_nu_OPCA_expansion}
\widehat\theta_\nu - \mbox{sign}\langle
\widehat\theta_\nu,\theta_\nu\rangle \theta_\nu = - H_\nu \mathbf{S}
\theta_\nu + R_\nu,
\end{equation}
where
\begin{equation}\label{eq:H_nu_Sigma}
H_\nu \equiv H_\nu(\Sigma) := \sum_{1\leq \mu \neq \nu \leq M}
\frac{1}{\lambda_\mu-\lambda_\nu} \theta_\mu \theta_\mu^T  -
\frac{1}{\lambda_\mu} P_\perp,
\end{equation}
where $P_\perp = I - \sum_{\mu=1}^M \theta_\mu \theta_\mu^T$. Note
that $H_\nu\theta_\nu = 0$ and that $H_\nu \Sigma \theta_\nu = 0$.
The key quantity in bounding the error term $R_\nu$ is
$$
\overline{\Delta}_\nu = \max\{(\lambda_\nu -
\lambda_{\nu+1})^{-1},(\lambda_{\nu-1}-\lambda_\nu)^{-1}\}
\parallel \mathbf{S} - \Sigma \parallel.
$$
Indeed, from (\ref{eq:eigenvec_error_leading}), when
$\overline{\Delta}_\nu < 1/4$, we have, for some constant $C > 0$,
$$
\parallel R_\nu \parallel \leq C \parallel H_\nu \mathbf{S} \theta_\nu \parallel \overline{\Delta}_\nu.
$$
Set $\delta_{n\nu}' = C \overline{\Delta}_\nu$. We will show that as
$n \to \infty$, $\delta_{n\nu}' \to 0$ with probability approaching
1 and
\begin{equation}\label{eq:L_theta_nu_hat_bounds}
\parallel H_\nu \mathbf{S}\theta_\nu \parallel^2 (1-\delta_{n\nu}')^2 \leq
L(\widehat\theta_\nu,\theta_\nu) \leq \parallel H_\nu \mathbf{S}\theta_\nu
\parallel^2 (1+\delta_{n\nu}')^2.
\end{equation}
Theorem \ref{thm:OPCA_risk_bound} then follows from an (exact,
non-asymptotic) evaluation
\begin{equation}\label{eq:expected_H_S_theta_nu}
\mathbb{E}\parallel H_\nu \mathbf{S}\theta_\nu \parallel^2  =
\frac{N-M}{nh(\lambda_\nu)} + \frac{1}{n} \sum_{\mu\neq \nu}
\frac{(1+\lambda_\mu)(1+\lambda_\nu)}{(\lambda_\mu -
\lambda_\nu)^2}~.
\end{equation}
We begin with the evaluation of (\ref{eq:expected_H_S_theta_nu}).
First we derive a convenient representation of $H_\nu
\mathbf{S}\theta_\nu $. In matrix form, model (\ref{eq:basic})
becomes
\begin{equation}\label{eq:basic_matrix_repr}
\mathbf{X} = \sum_{\mu=1}^M \sqrt{\lambda_\mu} \theta_\mu v_\mu^T + \mathbf{Z}.
\end{equation}
For $\nu=1,\ldots,M$, define
\begin{equation}\label{eq:z_w_nu_def}
z_\nu = \mathbf{Z}^T \theta_\nu, \qquad w_\nu = \mathbf{X}^T
\theta_\nu = \sqrt{\lambda_\nu}v_\nu + z_\nu.
\end{equation}
Define
\begin{equation}\label{eq:average_dot_product}
\langle \mathbf{a},\mathbf{b}\rangle_n := \frac{1}{n} \sum_{i=1}^n
a_ib_i ~~\mbox{for arbitrary} ~~\mathbf{a}, \mathbf{b} \in
\mathbb{R}^n.
\end{equation}
Then we have
$$
\mathbf{S}\theta_\nu = \frac{1}{n} \mathbf{X} w_\nu = \sum_{\mu=1}^M
\sqrt{\lambda_\mu} \langle v_\mu,w_\nu\rangle_n \theta_\mu + \frac{1}{n}
\mathbf{Z}w_\nu.
$$
Using (\ref{eq:z_w_nu_def}),
$$
H_\nu \mathbf{Z} w_\nu = \sum_{\mu\neq \nu} \frac{\langle z_\mu,
w_\nu\rangle}{\lambda_\mu - \lambda_\nu} \theta_\mu -
\frac{1}{\lambda_\nu} P_\perp \mathbf{Z}w_\nu.
$$
Using (\ref{eq:H_nu_Sigma}), $H_\nu \theta_\mu = (\lambda_\mu -
\lambda_\nu)^{-1}\theta_\mu$  for $\mu \neq \nu$, and we arrive at
the desired representation
\begin{equation}\label{eq:H_nu_Sigma_expansion}
H_\nu \mathbf{S}\theta_\nu = \sum_{\mu\neq \nu} \frac{\langle
w_\mu,w_\nu\rangle_n}{\lambda_\mu - \lambda_\nu} \theta_\mu -
\frac{1}{n\lambda_\nu} P_\perp \mathbf{Z}w_\nu.
\end{equation}
By orthogonality,
\begin{equation}\label{eq:H_nu_S_theta_nu_norm_squared}
\parallel H_\nu \mathbf{S}\theta_\nu \parallel^2 = \sum_{\mu\neq \nu} \frac{\langle w_\mu,w_\nu\rangle_n^2}
{(\lambda_\mu - \lambda_\nu)^2} + \frac{1}{n^2\lambda_\nu^2} w_\nu^T
\mathbf{Z}^T P_\perp \mathbf{Z}w_\nu.
\end{equation}
Now we compute the expectation. One verifies that $z_\nu \sim
N(0,I_n)$ independently of each other and of each $v_\nu \sim
N(0,I_n)$, so that $w_\nu \sim N(0,(1+\lambda_\nu) I_n)$
independently. Hence, for $\mu \neq \nu$,
\begin{eqnarray}\label{eq:expected_w_mu_w_nu}
\mathbb{E} \langle w_\mu,w_\nu \rangle_n^2 &=& n^{-2} \mathbb{E} \tr(w_\nu w_\nu^T w_\mu w_\mu^T) \nonumber\\
&=& n^{-2} \tr((1+\lambda_\mu)(1+\lambda_\nu) I_n) \nonumber\\
&=& n^{-1} (1+\lambda_\mu)(1+\lambda_\nu).
\end{eqnarray}
From (\ref{eq:z_w_nu_def}),
\begin{eqnarray*}
\mathbb{E}[w_\nu^T \mathbf{Z}^T P_\perp \mathbf{Z} w_\nu
|\mathbf{Z}] &=& z_\nu^T \mathbf{Z}^T P_\perp \mathbf{Z} z_\nu +
\lambda_\nu \mathbb{E}[v_\nu^T
\mathbf{Z}^T P_\perp \mathbf{Z} v_\nu|\mathbf{Z}] \\
&=& \tr(\mathbf{Z}\mathbf{Z}^T P_\perp \mathbf{Z}\mathbf{Z}^T
\theta_\mu\theta_\mu^T) + \lambda_\nu \tr(P_\perp
\mathbf{Z}\mathbf{Z}^T).
\end{eqnarray*}
Now, it can be easily verified that if $W := \mathbf{Z}\mathbf{Z}^T \sim$
$W_N(n,I)$, then for arbitrary symmetric $N\times N$ matrices $Q$, $R$, we
have,
\begin{equation}\label{eq:trace_quadratic}
\mathbb{E}\tr(WQWR) = n[\tr(QR) + \tr(Q)\tr(R)] + n^2 \tr(QR).
\end{equation}
Taking $Q = P_\perp$ and $R = \theta_\mu\theta_\mu^T$, by
(\ref{eq:trace_quadratic}) we have
\begin{equation}\label{eq:expected_w_nu_P_Z_quad}
\mathbb{E}[w_\nu^T \mathbf{Z} P_\perp \mathbf{Z}w_\nu] = n
\tr(P_\perp) + n \lambda_\nu \tr(P_\perp) = n(N - M)(1+\lambda_\nu).
\end{equation}
Combining (\ref{eq:expected_w_mu_w_nu}) with (\ref{eq:expected_w_nu_P_Z_quad})
in computing the expectation of (\ref{eq:H_nu_S_theta_nu_norm_squared}), we
obtain the expression (\ref{eq:expected_H_S_theta_nu}) for $\mathbb{E}\parallel
H_\nu \mathbf{S}\theta_\nu \parallel^2$.

\subsection*{Bound for $\parallel \mathbf{S} - \Sigma \parallel$}

We begin with the decomposition of the sample covariance matrix
$\mathbf{S}$. Introduce the abbreviation $\xi_\mu = n^{-1}
\mathbf{Z}v_\mu$. Then,
\begin{equation}\label{eq:S_expansion}
\mathbf{S} = \sum_{\mu=1}^M \sum_{\mu'=1}^M \sqrt{\lambda_\mu
\lambda_{\mu'}} \langle v_\mu,v_{\mu'}\rangle_n
\theta_\mu\theta_{\mu'}^T + \sum_{\mu=1}^M \sqrt{\lambda_\mu}
(\theta_\mu \xi_\mu^T + \xi_\mu \theta_\mu^T) + n^{-1}
\mathbf{Z}\mathbf{Z}^T
\end{equation}
and hence
\begin{eqnarray}\label{eq:S_Sigma_diff_prelim}
\parallel \mathbf{S} - \Sigma \parallel &\leq&
\sum_{\mu=1}^M \sum_{\mu'=1}^M \sqrt{\lambda_\mu \lambda_{\mu'}}
|\langle v_\mu, v_{\mu'}\rangle_n - \delta_{\mu\mu'}| \nonumber\\
&& + 2 \sum_{\mu=1}^M \sqrt{\lambda_\mu} \parallel \xi_\mu
\parallel +
\parallel n^{-1} \mathbf{Z}\mathbf{Z}^T - I \parallel,
\end{eqnarray}
where $\delta_{\mu\mu'}$ denotes the Kronecker symbol.
Let $D_1$ be the intersection of all the events (for some constant $c
> 0$):
\begin{eqnarray*}
D_{11} &:=&\{ |\parallel v_\mu \parallel_n^2 - 1| \leq 2c \sqrt{n^{-1} \log n}, ~1 \leq \mu \leq M\}, \\
D_{12} &:=& \{|\langle v_\mu, v_\nu \rangle_n | \leq c \sqrt{n^{-1} \log n}, ~1 \leq \mu \neq \mu' \leq M\}, \\
D_{13} &:=& \{\parallel \xi_\mu \parallel \leq (1+2c\sqrt{N^{-1} \log
n})\sqrt{\frac{N}{n}}, ~1 \leq \mu \leq M\}.
\end{eqnarray*}
Since $v_\nu \stackrel{i.i.d.}{\sim} N(0,I_n)$ independent of $\mathbf{Z}$, we
have $\mathbf{Z} v_\nu/\parallel v_\nu
\parallel \sim N(0,I_N)$ independently of $v_\nu$, and $\parallel v_\nu \parallel^2 \sim \chi_n^2$. Moreover,
$$
D_{11} \cap \{\parallel \mathbf{Z}v_\mu\parallel^2/\parallel v_\mu \parallel^2
\leq 1+2c\sqrt{N^{-1}\log n},~1\leq \mu \leq M\} \subset D_{13}.
$$
Hence, we use Lemmas \ref{lemma:chi_square_large_dev} and
\ref{lemma:quad_form_large_dev} to prove that
\begin{equation}\label{eq:D_1_prob_bound}
1-\mathbb{P}(D_1) \leq 3 M n^{-c^2} + M(M-1) n^{-(3/2)c^2 + O(n^{-1}
\log n)}.
\end{equation}
Define $D_2$ to be be the event that
\begin{equation}\label{eq:D_2_def}
D_2 := \left\{\parallel \frac{1}{n} \mathbf{Z}\mathbf{Z}^T - I_N
\parallel \leq \frac{N}{n} + 2 \sqrt{\frac{N}{n}} + ct_n\right\},
\end{equation}
with $t_n$ as in Lemma \ref{lemma:eigen_deviation_bound} with $N_n =
\max\{n,N\} = n$ so that $t_n = 8\sqrt{n^{-1} \log n}$. Lemma
\ref{lemma:eigen_deviation_bound} also establishes that $1 -
\mathbb{P}(D_2) \leq 2 n^{-c^2}$. Using the notation $\eta_n :=
(N^{-1} \log n)^{1/2}$, we have, on $D_1 \cap D_2$,
\begin{eqnarray}\label{eq:S_Sigma_diff_prob_bound}
\parallel \mathbf{S} - \Sigma \parallel &\leq& 2c(\sum_{\mu=1}^M \sqrt{\lambda_\mu})^2 \eta_n
+ 2 (\sum_{\mu=1}^M \lambda_\mu)(1+2c\eta_n) \sqrt{\frac{N}{n}} \nonumber\\
&& + 2\sqrt{\frac{N}{n}} + \frac{N}{n} + c t_n.
\end{eqnarray}
Recalling that $\rho_\nu = \lambda_\nu/\lambda_1$ for
$\nu=1,\ldots,M$, we have for large $n$ that
$$
\overline{\Delta}_\nu \leq C_\nu(\rho) \frac{\parallel \mathbf{S} -
\Sigma \parallel}{\lambda_1},
$$
where, say $C_\nu(\rho) = 2 \max\{(\rho_\nu -
\rho_{\nu+1})^{-1},(\rho_{\nu-1}-\rho_\nu)^{-1}\}$. Observe that
$t_n/\lambda_1 = 8\eta_n \sqrt{N/(n\lambda_1)^2}$. Now, substitute
(\ref{eq:S_Sigma_diff_prob_bound}) to conclude that there are
functions $B_i(\rho)$ such that on $D_n := D_1 \cap D_2$,
$$
\overline{\Delta}_\nu  \leq B_1(\rho) \eta_n + B_2(\rho)
(1+2c\eta_n) \sqrt{\frac{N}{n\lambda_1}} +
2\sqrt{\frac{N}{n\lambda_1^2}} + \frac{N}{n\lambda_1} +
8c\eta_n\sqrt{\frac{N}{n\lambda_1^2}}~.
$$
Our assumptions imply that
$$
\eta_n = \sqrt{\frac{\log n}{N}} \to 0 \qquad \mbox{and} \qquad
\frac{N}{n\lambda_1^2} + \frac{N}{n\lambda_1} = \frac{N}{nh(\lambda_1)} \to 0,
$$
so that $\overline{\Delta}_\nu \to 0$. To summarize, choose $c =
\sqrt{2}$, say, so that on $D_n$, which has probability at least
$1-O(n^{-2})$, we have $\delta_{n\nu}'\to 0$. This completes the
proof of (\ref{eq:L_theta_nu_hat_bounds}).

Theorem \ref{thm:OPCA_risk_bound} now follows from noticing that
$L(\widehat\theta_\nu,\theta_\nu) \leq 2$ and so
\begin{equation*}
\mathbb{E}[L(\widehat\theta_\nu,\theta_\nu), (D_1 \cap D_2)^c] \leq 2
\mathbb{P}((D_1 \cap D_2)^c) = O(N_n^{-2}) = o(\mathbb{E}\parallel H_\nu
\mathbf{S}\theta_\nu
\parallel^2),
\end{equation*}
and an additional computation using
(\ref{eq:H_nu_S_theta_nu_norm_squared}) which shows that
\begin{equation*}
\mathbb{E}[\parallel H_\nu \mathbf{S} \theta_\nu \parallel^2, D_n^c]
\leq (\mathbb{E}[\parallel H_\nu \mathbf{S} \theta_\nu
\parallel^4)^{1/2} P(D_n^c) = o(\mathbb{E}[\parallel H_\nu \mathbf{S} \theta_\nu
\parallel^2).
\end{equation*}

\subsection{Lower bound on the minimax risk} \label{sec:highdpca_proof_minimax_lbd}

In this subsection, we prove Theorems \ref{th:three-way-lower} and
\ref{th:sparse-lower}.  The key idea in the proofs is to utilize the
geometry of the parameter space in order to construct appropriate
finite dimensional subproblems for which bounds are easier to
obtain. We first give an overview of the general machinery used in
the proof.

\subsubsection*{Risk bounding strategy}\label{subsec:highdpca_lower_bound_strategy}

A key tool for deriving lower bounds on the minimax risk is
\textit{Fano's Lemma}. In this subsection, we use superscripts on
vectors $\theta$ as indices, not exponents. First,
we construct a large finite subset ${\cal F}$ of $\Theta_q^M(C_1,\ldots,C_M)$, such that the
following property holds, for a given $\nu \in \{1,\ldots,M\}$.
\begin{itemize}
\item[]
If $\boldsymbol\theta^{1}, \boldsymbol\theta^{2} \in {\cal F}$, then
$L(\theta_\nu^{1},\theta_\nu^{2}) \geq 4\delta$, for some
    $\delta > 0$ (to be chosen).
\end{itemize}
This property will be referred to as ``$4\delta$-distinguishability in
$\theta_\nu$''. Given any estimator $\widehat{\boldsymbol\theta}$ of
$\boldsymbol\theta$, based on data $\mathbf{X}_n = (X_1,\ldots,X_n)$, define a
new estimator
$\phi(\mathbf{X}_n) = \bs \theta^*$, whose $M$ components are given by
$\theta^*_\nu = \arg \min_{\boldsymbol\theta \in {\cal F}} L(\widehat
\theta_\nu,\theta_\nu)$, where $\widehat \theta_\nu$ is the $\nu$-th column of
$\widehat{\boldsymbol\theta}$. Then, by Chebyshev's inequality and the
$4\delta$-distinguishability in $\theta_\nu$, it follows that
\begin{eqnarray}
            \label{eq:risk_lb_M}
    \sup_{\boldsymbol\theta \in \Theta_q^M(C_1,\ldots,C_M)}
\mathbb{E}_{\boldsymbol\theta} L(\widehat\theta_\nu,\theta_\nu)
&\geq& \delta \sup_{\boldsymbol\theta \in {\cal F}}
\mathbb{P}_{\boldsymbol\theta} ( \phi(\mathbf{X}_n) \neq
\boldsymbol\theta).
\end{eqnarray}
The task is then to find an appropriate lower bound for the quantity on the
right hand side of (\ref{eq:risk_lb_M}). For this, we use the following version
of Fano's lemma, due to \cite{Birge2001}, modifying a result of
\cite{YangB1999} (p. 1570-71).

\begin{lemma} \label{lemma:Fano_lb_large}
Let $\{P_\theta : \theta \in \Theta\}$ be a family  of probability
distributions on a common measurable space, where $\Theta$ is an arbitrary
parameter set.
Let $p_{max}$ be the minimax risk over $\Theta,$ with the loss function $L'(\theta,\theta') = \mathbf{1}_{\theta \neq \theta'}$,

$$
p_{max} = \inf_T  \sup_{\theta \in \Theta} \mathbb{P}_\theta(T \neq \theta) =
\inf_T  \sup_{\theta \in \Theta} \mathbb{E} L'(\theta,T),
$$
where $T$ denotes an arbitrary estimator of $\theta$ with values in $\Theta$.
Then for any finite subset ${\cal F}$ of $\Theta$, with elements
$\theta_1,\ldots,\theta_J$ where $J = |{\cal F}|$,
\begin{equation}\label{eq:Fano}
p_{max} \geq 1 -  \inf_{Q} ~\frac{J^{-1} \sum_{i=1}^J K(P_i, Q) + \log  2}{\log
J}
\end{equation}
where $P_i =  \mathbb{P}_{\theta_i}$, and $Q$ is an arbitrary probability
distribution, and $K(P_i,Q)$ is the Kullback-Leibler divergence of $Q$ from
$P_i$.
\end{lemma}

The following lemma, proven in Section \ref{sec:auxiliary_results}, gives the Kullback-Leibler discrepancy corresponding to
two different values of the parameter.

\begin{lemma}\label{prop:multi_KL_div}
Let $\boldsymbol{\theta}^{j} := [\theta_1^{j} : \ldots :
\theta_M^{j}]$, $j=1,2$ be two parameters (i.e., for each $j$,
$\theta_k^j$'s are orthonormal). Let $\Sigma_{j}$ denote the matrix
given by (\ref{eq:Sigma_basic}) with $\bs\theta = \bs\theta^{j}$
(and $\sigma = 1$). Let $P_j$ denote the joint probability
distribution of $n$ i.i.d. observations from $N(0,\Sigma_{j})$. Then
the Kullback-Leibler discrepancy of $P_2$ with respect to $P_1$ is
given by
\begin{equation}\label{eq:multi_KL_div}
{\cal K}_{1,2} := K(\boldsymbol\theta^{1},\boldsymbol\theta^{2}) =
\frac{n}{2} \Bigg[\sum_{\nu=1}^M \eta(\lambda_\nu) \lambda_\nu -
\sum_{\nu=1}^M\sum_{\mu=1}^M  \eta(\lambda_\nu) \lambda_{\mu}
|\langle \theta_{\mu}^{1}, \theta_\nu^{2}\rangle|^2\Bigg],
\end{equation}
where $\eta(\lambda) = \lambda/(1+\lambda)$.
\end{lemma}

\subsubsection*{Geometry of the hypothesis set and Sphere Packing}\label{subsec:highdpca-geometry}

Next, we describe the construction of a large set of hypotheses ${\cal F}$, satisfying the $4\delta$
distinguishability condition. Our construction is based on the well studied sphere packing problem,
namely how many unit vectors can be packed onto $\mathbb{S}^{m-1}$, with given minimal pairwise distance
between any two vectors.

Here we follow the construction due to
\cite{Zong1999} (p. 77). Let $m$ be a large
positive integer, and $m_0 = \lfloor 2m/9\rfloor$. Define $Y_m^*$ as
the maximal set of points of the form ${\bf z}=(z_1,\ldots,z_m)$ in
$\mathbb{S}^{m-1}$ such that the following is true:
\begin{equation*}\label{eq:separation}
\sqrt{m_0} z_i \in \{-1,0,1\} ~\forall~ i, ~~ \sum_{i=1}^m |z_i| =
\sqrt{m_0} ~~\mbox{and, for} ~~{\bf z}, {\bf z}' \in Y_m^*, ~~
\parallel {\bf z} - {\bf z}'\parallel \geq 1.
\end{equation*}
For any $m \geq 1$, the maximal number of points lying on $\mathbb{S}^{m-1}$
such that any two points are at distance at least 1, is called the
\textit{kissing number} of an $m$-sphere.
\cite{Zong1999} uses the
construction described above to derive a lower bound on the
\textit{kissing number}, by showing that $|Y_m^*| \geq
(9/8)^{m(1+o(1))}$ for $m$ large.

Next, for $m<N-M$ we use the sets $Y_m^*$ to construct our hypothesis set $\cal F$
of same size, $|{\cal F}| = |Y_m^*|$. To this end, let $\{{\bf e}_\mu\}_{\mu=1}^N$
denote the standard basis of $\mathbb{R}^N$. Our initial set ${\bs \theta}^0$ is composed
of the first $M$ standard basis vectors, ${\bs \theta}^0 = [{\bf e}_1:\ldots:{\bf e}_M]$.
Then, for fixed $\nu$, and values of $m,r$ yet to be determined, each of the other hypotheses ${\bs \theta}^j\in{\cal F}$ has the same vectors as ${\bf \theta}^0$ for $k\neq \nu$.
The difference is that the $\nu$-th vector is instead given by
\begin{equation}
                \label{eq:theta_nu_j_F}
    \theta_\nu^{j} = \sqrt{1-r^2} ~\mathbf{e}_\nu + r \sum_{l=1}^m
        z_l^{j} \mathbf{e}_{M+l}, ~~~j=1,\ldots,|{\cal F}|,
\end{equation}
where $\mathbf{z}^{j}= (z_1^{j},\ldots,z_m^{j})$, $j\geq 1$, is an
enumeration of the elements of $Y_m^*$. Thus $\theta_\nu^{j}$
perturbs $\mathbf{e}_\nu$ in subsets of the fixed set of coordinates $\{M+1, \dots, M+m\}$, according to the sphere
packing construction for $\mathbb{S}^{m-1}$.

The
construction ensures that $\theta_1^j,\ldots,\theta_M^j$ are
orthonormal for each $j$.
Furthermore, (\ref{eq:multi_KL_div}) simplifies to
\begin{equation}\label{eq:KL_div}
K(\boldsymbol\theta^{j},\boldsymbol\theta^{0}) = \frac{1}{2}
nh(\lambda_\nu) (1 - (\langle
\theta_\nu^{j},\theta_\nu^{0}\rangle)^2) = \frac{1}{2}
nh(\lambda_\nu) r^2, ~~j=1,\ldots,|{\cal F}|.
\end{equation}
Finally, by construction, for any $\bs \theta^j,\bs \theta^k\in {\cal F}$ with $j\neq k$
\begin{equation}\label{eq:r_square_distinguish}
L(\theta_\nu^{j},\theta_\nu^{k}) \geq r^2,
\end{equation}
In other words, the set ${\cal F}$ is $r^2$-distinguishable in $\theta_\nu$.
Consequently, combining (\ref{eq:risk_lb_M}) and (\ref{eq:KL_div}),
\begin{equation}
  \label{eq:rstarbd}
  R_\nu^* = \inf_{\hat \theta_\nu} \sup_{\Theta_q(\mathbf{C})}
  \mathbb{E} L(\hat \theta_\nu, \theta_\nu)
   \geq (r^2/4) [1 -  a(r, \mathcal{F})],
\end{equation}
with
\begin{equation}\label{eq:a_r_F_0}
  a(r,\mathcal{F}) = \frac{\tfrac{1}{2} n h(\lambda_\nu) r^2 + \log 2}{\log
    |\mathcal{F}|}~.
\end{equation}



\subsubsection*{Proof of  Theorem \ref{th:three-way-lower}}


Let $m$ be an integer yet to be specified and let $r \in (0,1)$. Let $Y_m^*$ be
the sphere-packing set defined above, and let $\cal F$ be the corresponding set
of hypotheses, defined via (\ref{eq:theta_nu_j_F}).

Let $c_1 = \log(9/8)$, then we have $\log |\mathcal{F}| \geq b_m
c_1 m$, where $b_m \to 1$ as $m\to\infty$.
Inserting the following value for $r=r(m)$,
\begin{equation}
  \label{eq:rdef}
  r^2 = \frac{c_1 m}{n h(\lambda_\nu)},
\end{equation}
into Eq. (\ref{eq:a_r_F_0}) gives that
\begin{displaymath}
  a(r, \mathcal{F}) \leq
    \frac{ \tfrac{1}{2} c_1 m + \log 2}{b_m c_1 m}~.
\end{displaymath}
Therefore, so long as $m \geq m_*$, an absolute constant, we have
$a(r, \mathcal{F}_0) \leq 3/4$.

We need to ensure that $\theta_\nu^{j} \in \Theta_q(C_\nu)$. Since
exactly $m_0$ coordinates are non-zero out of $\{ M+1, \dots,
M+m\}$,
\begin{displaymath}
  \| \theta_\nu^{j} \|_q^q = (1 - r^2)^{q/2} + r^q m_0^{1-q/2} \leq
  1 + a_q r^q m^{1-q/2}
\end{displaymath}
where $a_q = (2/9)^{1-q/2}$. A sufficient condition for
$\theta_\nu^{(j)} \in \Theta_q(C_\nu)$ is that
\begin{equation}
  \label{eq:int-cond}
  a_q m (r^2/m)^{q/2} \leq \bar C_\nu^q.
\end{equation}
Substituting (\ref{eq:rdef}) puts this into the form
\begin{displaymath}
   m \leq \frac{1}{a_q c_1^{q/2}} \bar C_\nu^q [n h(\lambda_\nu)]^{q/2}.
\end{displaymath}

To simultaneously ensure that (i) $r^2 < 1$, (ii) $m$ does not
exceed the number of available co-ordinates, $N-M$, and (iii)
$\theta_\nu^{j} \in \Theta_q(C_\nu)$, we set
\begin{displaymath}
  m = \min \{ \lfloor n h(\lambda_\nu)\rfloor, N-M, \lfloor A_q \bar C_\nu^q (n
  h(\lambda_\nu))^{q/2}\rfloor \},
\end{displaymath}
where $A_q = 1/(a_q c_1^{q/2})$. Recalling the notations
(\ref{eq:tau-nu2}), (\ref{eq:m-nu}) and (\ref{eq:Nprime}), this
becomes (without loss of generality assuming $nh(\lambda_\nu)$ and
$m_\nu$ to be integers)
\begin{displaymath}
  m = \min \{ \tau_\nu^{-2}, N', m_\nu \} = \tau_\nu^{-2} \min\{1,\tau_\nu^2 \cdot \min\{N', m_\nu\} \}
\end{displaymath}
and Theorem \ref{th:three-way-lower} follows.

\subsubsection*{Proof of Theorem \ref{th:sparse-lower}}

The construction of the set of hypotheses in the proof of
Theorem \ref{th:three-way-lower} considered a fixed set of potential
non-zero coordinates, namely $\{M+1,\ldots,M+m\}$. However, in the
\textit{ultra-sparse} setting, when the effective
dimension is significantly smaller than the nominal dimension $N$, it is
possible to construct a much larger collection of hypotheses by allowing the set
of non-zero coordinates to span all remaining coordinates $\{M+1,\ldots,N\}$.

In the proof of Theorem \ref{th:sparse-lower} we shall use the following lemma,
proven in Section \ref{sec:auxiliary_results}. Call $A \subset \{ 1, \ldots, N
\}$ an \textit{$m-$set} if $|A| = m$.

\begin{lemma} \label{lemma:counting}
  Let $k$ be fixed, and let $\mathcal{A}_k$ be the maximal collection of
  $m-$sets such that the intersection of any two members has
  cardinality at most $k-1$. Then, necessarily,
  \begin{equation}
    \label{eq:cardinality}
    | \mathcal{A}_k | \geq \binom{N}{k} \bigg/ \binom{m}{k}^2.
  \end{equation}
Let $k = [m_0/2] + 1$ and $m_0 = [\beta m]$ with $0 < \beta < 1.$
Suppose that $m, N \rightarrow \infty$ with $m = o(N)$. Then
\begin{equation}
  \label{eq:entropy}
  | \mathcal{A}_k | \geq \exp [ N \mathcal{E}(\beta m/2N) - 2 m
  \mathcal{E}(\beta/2) ] (1 + o(1)).
\end{equation}
where ${\cal E}(x)$ is the Shannon entropy function,
$$
{\cal E}(x) = - x \log(x) - (1-x) \log(1-x), ~~ 0 < x < 1.
$$
\end{lemma}

Let $\pi$ be an $m-$set contained in $\{ M+1, \dots, N \}$, and
construct a family $\mathcal{F}_\pi$ by modifying
\eqref{eq:theta_nu_j_F} to use the set $\pi$ rather than the fixed
set $\{M+1, \dots, M+m\}$ as in Theorem \ref{th:three-way-lower}:
\begin{equation*}\label{eq:theta_nu_j_F_pi}
\theta_\nu^{(j,\pi)} = \sqrt{1-r^2} ~\mathbf{e}_\nu + r \sum_{l \in
\pi} z_l^{j} \mathbf{e}_l, ~~~j=1,\ldots,|Y_m^*|.
\end{equation*}
We will choose $m$ below to ensure that $\theta_\nu^{(j,\pi)} \in
\Theta_q(C_\nu)$.
Let ${\cal P}$ be a collection of sets $\pi$ such that, for any
two sets $\pi$ and $\pi'$ in ${\cal P}$, the set $\pi \cap \pi'$ has
cardinality at most $m_0/2$. This ensures that the sets ${\cal
F}_\pi$ are disjoint  for $\pi \neq \pi'$, since each
$\theta_\nu^{(j,\pi)}$
is nonzero in exactly $m_0+1$ coordinates. This construction also
ensures that
$$
\mbox{for all} ~~ \mathbf{y}, \mathbf{y}' \in \bigcup_{\pi \in \cal
P} {\cal F}_\pi, \quad L(\mathbf{y},\mathbf{y}') \geq
\left(\frac{m_0}{2} +
\frac{m_0}{2}\right)\left(\frac{r}{\sqrt{m_0}}\right)^2 = r^2.
$$
Define ${\cal F} := \bigcup_{\pi \in \cal P} {\cal F}_\pi$. Then
\begin{equation}\label{eq:big_F_0_lower}
|{\cal F}| = |\bigcup_{\pi \in \cal P} {\cal F}_\pi| = |{\cal P}|
~|Y_m^*| \geq |{\cal P}| (9/8)^{m(1+o(1))}.
\end{equation}
By Lemma \ref{lemma:counting}, there is a collection ${\cal P}$ such
that $|{\cal P}|$ is at least $\exp( [N {\cal E}(m/9N) - 2m{\cal
E}(1/9)](1+o(1)))$.
Since ${\cal E}(x) \geq -x\log x$, it follows from
(\ref{eq:big_F_0_lower}) that,
$$
\frac{\log |{\cal F}|}{m} \geq  \left( \frac{1}{9} \log
\frac{9N}{m} - 2 {\cal
  E}(1/9)\right) + \log(9/8)(1+o(1))
    \geq \frac{\alpha}{9} \log N  +O(1),
$$
since $m = O(N^{1-\alpha})$.

Proceeding as for Theorem \ref{th:three-way-lower}, we have $\log
|\mathcal{F}| \geq b_m (\alpha/9)  m \log N$, where $b_m \to 1$.
Let us set (with $m$ still to be specified)
\begin{equation}
  \label{eq:rdef_refined}
  r^2 = m \frac{(\alpha /9) \log N}{n h(\lambda_\nu)}
      = m \bar \tau_\nu^2,
\end{equation}

Again, we need to ensure that $\theta_\nu^{(j,\pi)} \in
\Theta_q(C_\nu)$, which as before is implied by \eqref{eq:int-cond}.
Substituting (\ref{eq:rdef_refined}) puts this into the form
\begin{displaymath}
   m \leq \bar{m}_\nu = a_q^{-1} (\bar C_\nu/ \bar \tau_\nu)^q.
\end{displaymath}
To simultaneously ensure that (i) $r^2 < 1$; (ii) $m$ does not
exceed the number of available co-ordinates, $N-M$; and (iii)
$\theta_\nu^{j} \in \Theta_q(C_\nu)$, we set
\begin{displaymath}
m = \min \{ \lfloor \bar \tau_\nu^{-2}\rfloor, N-M, \lfloor a_q^{-1}
(\bar C_\nu^q/ \bar \tau_\nu)^q\rfloor \}.
\end{displaymath}
As $n,N\to\infty$, we have that $m = \lfloor a_q^{-1}
(\bar C_\nu/ \bar \tau_\nu)^q\rfloor$, and Theorem
\ref{th:sparse-lower} follows.

\subsection{Lower bound on the risk of the D.T. estimator}

To prove Theorem \ref{thm:diagonal_thresholding_risk}, assume w.l.g. that
$\langle \widehat \theta_{1,DT},\theta_1 \rangle
> 0$, and decompose the loss as
\begin{equation}\label{eq:loss_decomp}
L(\widehat \theta_{1,DT},\theta_1) = \parallel \theta_1 -
\theta_{1,I}\parallel^2 + \parallel \widehat \theta_{1,DT} -
\theta_{1,I}\parallel^2,
\end{equation}
where $I=I(\gamma_n)$ is the set of coordinates selected by the D.T.
scheme and $\theta_{1,I}$ denotes the subvector of $\theta_1$
corresponding to this set. Note that, in (\ref{eq:loss_decomp}), the
first term on the right can be viewed as a bias term while the
second term can be seen as a variance term.

We choose a particular vector $\theta_1 = \theta_{*} \in
\Theta_q(C)$ so that
\begin{equation}\label{eq:diagonal_thresh_bias}
\mathbb{E}\parallel \theta_* - \theta_{*,I}\parallel^2 \geq K
\bar{C}^q n^{-\frac{1}{2}(1-q/2)}.
\end{equation}
This, together with (\ref{eq:loss_decomp}), proves Theorem
\ref{thm:diagonal_thresholding_risk} since the worst case risk is
clearly at least as large as (\ref{eq:diagonal_thresh_bias}).
Accordingly, set $r_n = \bar{C}^{q/2} n^{-\frac{1}{4}(1-q/2)}$,
where $\bar{C}^q = C^q -1$. Since $C^q n^{q/4} = o(n^{1/2})$, we
have $r_n = o(1)$, and so for sufficiently large $n$, we can take
$r_n < 1$ and define
\begin{eqnarray*}
\theta_{*,k} = \begin{cases} \sqrt{1-r_n^2} & ~\mbox{if}~~ k=1\\
\frac{r_n}{\sqrt{m_n}} &~\mbox{if}~~ 2
\leq k \leq m_n+1\\
0 & ~\mbox{if}~~ m_n+2 \leq k \leq N
\end{cases}
\end{eqnarray*}
where $m_n = \lfloor(1/2) \bar{C}^q n^{q/4}\rfloor$. Then by construction  $\theta_* \in \Theta_q(C)$, since
$$
\sum_{k=1}^N |\theta_{*,k}|^q = (1-r_n^2)^{q/2} + r_n^q m_n^{1-q/2}
< 1 + r_n^q m_n^{1-q/2} \leq 1 + \frac{\bar{C}^q}{2^{1-q/2}} < C^q,
$$
where the last inequality is due to $q \in (0,2)$ and $\bar{C}^q =
C^q -1$.

For notational convenience, let $\alpha_n = \gamma\sqrt{\log N/n}$.
Recall that D.T. selects all coordinates $k$ for which
$\mathbf{S}_{kk} > 1+\alpha_n$. Therefore, coordinate $k$ is
\textit{not selected} with probability
\begin{equation}\label{eq:DT_selection_prob}
p_k = \mathbb{P}(\mathbf{S}_{kk} < 1+\alpha_n) =
\mathbb{P}\left(\frac{W_n}{n} < \frac{1+\alpha_n}{1+\lambda_1
\theta_{*,k}^2}\right)
\end{equation}
where $W_n \sim \chi_{n}^2$. Notice that, for $k = 2,\ldots,m_n+1$,
$p_k = p_2$, and $\theta_{*,k} =0$ for $k > m_n+1$. Hence,
\begin{equation*}
\mathbb{E}\parallel \theta_* - \theta_{*,I}\parallel^2 =
\sum_{k=1}^N p_k |\theta_{*,k}|^2 > p_2 \sum_{k=2}^{m_n+1}
|\theta_{*,k}|^2 = p_2 r_n^2 = p_2 \bar{C}^q
n^{-\frac{1}{2}(1-q/2)}.
\end{equation*}
Thus, to finish the proof of Theorem
\ref{thm:diagonal_thresholding_risk}, it is enough to show that $p_2
> 1-A_n$ for some $A_n$ that converges to 0 as $n \to \infty$.
Rewrite (\ref{eq:DT_selection_prob}) as
\begin{equation*}
p_k = \mathbb{P}\left(\frac{W_n}{n} < 1+\epsilon_k\right) = 1-
\mathbb{P}\left(\frac{W_n}{n} \geq 1+\epsilon_k\right)
~~\mbox{where}~~\epsilon_k =
\frac{\alpha_n-\lambda_1|\theta_{*,k}|^2}{1+\lambda_1
|\theta_{*,k}|^2}~.
\end{equation*}
Since $|\theta_{*,2}|^2 = r_n^2/m_n = 2n^{-1/2}(1+o(1))$, it follows
that
\begin{eqnarray*}
\epsilon_2 = \frac{\gamma \sqrt{\frac{\log N}{n}} - \lambda_1
\frac{r_n^2}{m_n}}{1+\lambda_1 \frac{r_n^2}{m_n}} =
\frac{1}{\sqrt{n}} \left(\frac{\gamma \sqrt{\log N} -
2\lambda_1}{1+2\lambda_1/\sqrt{n}}\right)(1+o(1))
\end{eqnarray*}
so that $n \epsilon_2^2 \to \infty$ as $n \to \infty$. This,
together with (\ref{eq:large_dev_chisq_2}), shows that $p_2 \geq
1-A_n$ where we can choose $A_n = \exp(-3n\epsilon_2^2/16) = o(1)$.

\section{Proof of relevant lemmas}\label{sec:auxiliary_results}

\subsection{Proof of Lemma \ref{lemma:eigen_deviation_bound}}

We use the following result on extreme eigenvalues of Wishart matrices by
\cite{Davidson2001}.
\begin{lemma}\label{lemma:extreme_singval_concen}
Let $Z$ be a $p \times q$ matrix of i.i.d. $N(0,1)$ entries with $p \leq q$.
Let $s_{max}(Z)$ and $s_{min}(Z)$  denote the largest and the smallest singular
value of $Z$, respectively. Then,
\begin{eqnarray}
\mathbb{P}(s_{max}(\frac{1}{\sqrt{q}} Z) > 1+ \sqrt{p/q} + t) &\leq&
e^{-qt^2/2} ,
\label{eq:extreme_singular_1}\\
\mathbb{P}(s_{min}(\frac{1}{\sqrt{q}}Z) < 1- \sqrt{p/q} - t) &\leq&
e^{-qt^2/2}. \label{eq:extreme_singular_2}
\end{eqnarray}
\end{lemma}

We apply Lemma \ref{lemma:extreme_singval_concen} separately for $N \leq n$ and
for $N > n$. Observe first that,
$$
\Delta := \parallel \frac{1}{n} \mathbf{Z}\mathbf{Z}^T - I_N
\parallel  = \max\{\lambda_1(n^{-1}\mathbf{Z}\mathbf{Z}^T) -1, 1-
\lambda_N(\mathbf{Z}\mathbf{Z}^T)\}.
$$
Consider first $N \leq n$ and let $s_\pm$ denote the maximum and minimum
singular values of $n^{-1/2} \mathbf{Z}$. Define $\gamma(t) := \sqrt{N/n} + t$
for $t > 0$. Then, since $\Delta = \max\{s_+^2 - 1, 1-s_-^2\}$,  and letting
$\Delta_n(t) := 2\gamma(t) + \gamma(t)^2$ we have
$$
\{\Delta > \Delta_n(t)\} \subset \{s_+ > 1+\gamma(t)\} \cup \{s_- <
1-\gamma(t)\}.
$$
Now, applying Lemma \ref{lemma:extreme_singval_concen} with $p=N$ and $q=n$, we
get
$$
\mathbb{P}( \Delta > \Delta_n(t)) \leq 2e^{-nt^2/2}.
$$
We observe that
\begin{equation}\label{eq:Delta_n_t}
\Delta_n(t) = (N/n + 2\sqrt{N/n}) + t(2+t+2\sqrt{N/n}).
\end{equation}
Now consider $N > n$. Noting that $\lambda_N(n^{-1} \mathbf{Z}\mathbf{Z}^T) =
0$, we have
$$
\Delta = \max\{(N/n)s_+^2 -1, 1\}.
$$
This time, let $\bar\gamma(t) := \sqrt{n/N} + t$ and $\Delta_N(t) :=
\max\{(N/n) (1+\bar\gamma(t))^2 - 1, 1\}$. We apply Lemma
\ref{lemma:extreme_singval_concen} with $p=n$, $q=N$, so that
$$
\mathbb{P}( \Delta > \Delta_N(t)) = \mathbb{P}(s_+ > 1+\bar\gamma(t)) \leq
e^{-nt^2/2},
$$
and observe that
\begin{equation}\label{eq:Delta_N_t}
\Delta_N(t) = (N/n + 2\sqrt{N/n}) + (N/n) t(2+t+2\sqrt{n/N}).
\end{equation}
Thus from (\ref{eq:Delta_n_t}) and (\ref{eq:Delta_N_t}), we have
$$
\Delta_{\max\{n,N\}}(t) \leq (N/n + 2\sqrt{N/n})  + t(N_n/n)(4+t).
$$
Now choose $t = c\sqrt{2\log N_n/N_n}$ so that tail probability is at most $2
e^{-N_n^2 t^2/2} = 2 N_n^{-c^2}$. The result is now proved, since if $c
\sqrt{\log n/n} \leq 1$ then $t(N_n/n)(4+t) \leq ct_n$.


\subsection{Proof of Lemma \ref{prop:multi_KL_div}}

Recall that, if distributions $F_1$ and $F_2$ have density functions $f_1$ and
$f_2$, respectively, such that the support of $f_1$ is contained in the support
of $f_2$, then the Kullback-Leibler discrepancy of $F_2$ with respect to $F_1$,
to be denoted by $K(F_1,F_2)$, is given by
\begin{equation}\label{eq:KL_general}
K(F_1,F_2) = \int \log \frac{f_1(y)}{f_2(y)} f_1(y) dy.
\end{equation}
For $n$ i.i.d. observations $X_i, i=1,\ldots,n$, the Kullback-Leibler
discrepancy is just $n$ times the Kullback-Leibler discrepancy for a single
observation. Therefore, without loss of generality we take $n = 1$. Since
\begin{equation}\label{eq:Sigma_inverse}
\Sigma^{-1} = (I - \sum_{\nu=1}^M \eta(\lambda_\nu) \theta_\nu \theta_\nu^T),
\end{equation}
the log-likelihood function for a single observation is given by
\begin{eqnarray}\label{eq:multi_log_likelihood}
\log f(x|\bs\theta) &=& -\frac{N}{2}\log(2\pi) - \frac{1}{2}\log|\Sigma|
- \frac{1}{2} x^T\Sigma^{-1}x \nonumber\\
&=& -\frac{N}{2}\log(2\pi) - \frac{1}{2}\sum_{\nu = 1}^M
\log(1+\lambda_\nu)  \nonumber\\
&& - \frac{1}{2} \left(\langle x,x \rangle - \sum_{\nu=1}^M \eta(\lambda_\nu)
\langle x,\theta_\nu \rangle^2\right).
\end{eqnarray}
From (\ref{eq:multi_log_likelihood}), we have
\begin{eqnarray*}
&& {\cal K}_{1,2} \\
&=& \mathbb{E}_{\bs\theta^{1}} \left(\log f(X|\bs\theta^{1}) - \log
f(X|\bs\theta^{2})\right) \nonumber\\
&=& \frac{1}{2}  \sum_{\nu=1}^M \eta(\lambda_\nu)
[\mathbb{E}_{\bs\theta^{1}}(\langle X , \theta_\nu^{1} \rangle)^2 -
\mathbb{E}_{\bs\theta^{1}}(\langle X ,
\theta_\nu^{2} \rangle)^2] \nonumber\\
&=& \frac{1}{2}  \sum_{\nu=1}^M \eta(\lambda_\nu) [\langle \theta_\nu^{1},
\Sigma_{(1)} \theta_\nu^{1}\rangle - \langle
\theta_\nu^{2}, \Sigma_{(1)} \theta_\nu^{2}\rangle] \nonumber\\
&=& \frac{1}{2}  \sum_{\nu=1}^M \eta(\lambda_\nu) \left[(\parallel\theta_\nu^1
\parallel^2 - \parallel\theta_\nu^2
\parallel^2)
 + \sum_{\mu=1}^M \lambda_{\mu} \{(\langle
\theta_{\mu}^{1},\theta_\nu^{1}\rangle)^2 - (\langle \theta_{\mu}^{1},
\theta_\nu^{2}\rangle)^2\}\right],
\end{eqnarray*}
which equals the RHS of (\ref{eq:multi_KL_div}), since the columns of
$\bs\theta^{j}$ are orthonormal for each $j=1,2$.


\subsection{Proof of Lemma \ref{lemma:counting}}

Let $\mathcal{P}_m$ be the collection of all $m-$sets of $\{ 1,  \ldots,
  N \}$, clearly $|\mathcal{P}_m| = \binom{N}{m}.$
  For any $m-$set $A$, let $\mathcal{I}(A)$ denote the collection of
  ``inadmissible''
  $m-$sets $A'$ for which $| A \cap A' | \geq k$. Clearly
  \begin{displaymath}
    | \mathcal{I}(A) | \leq \binom{m}{k} \binom{N-k}{m-k}.
  \end{displaymath}
If $\mathcal{A}_k$ is maximal, then $\mathcal{P}_m = \cup_{A \in
  \mathcal{A}_k} \mathcal{I}(A)$, and so \eqref{eq:cardinality}
follows from the inequality
\begin{displaymath}
  |\mathcal{P}_m| \leq | \mathcal{A}_k| \, \max_A |\mathcal{I}(A)|,
\end{displaymath}
and rearrangement of factorials.

Turning to the second part, we recall that Stirling's formula shows that if $k$
and $N \rightarrow \infty$,
\begin{displaymath}
  \binom{N}{k} = \theta \bigg( \frac{N}{2 \pi k (N-k)} \bigg)^{1/2}
                    \exp \Big\{ N \mathcal{E} \Big(\frac{k}{N} \Big) \Big\},
\end{displaymath}
where $\theta \in ( 1 - (6k)^{-1}, 1 + (12 N)^{-1})$. The coefficient
multiplying the exponent in $ \binom{N}{k} \big/ \binom{m}{k}^2$ is
\begin{displaymath}
  \sqrt{2 \pi k} (1 - k/N)^{-1/2} (1 - k/m)
     \sim \sqrt{ \pi \beta m} ( 1- \beta/2) \rightarrow \infty
\end{displaymath}
under our assumptions, and this yields \eqref{eq:entropy}.


\end{document}